\documentclass[a4paper, 11pt]{article}

\usepackage{srcltx}
\usepackage{indentfirst}
\usepackage{epsfig}
\usepackage{psfrag}
\usepackage{graphicx}
\usepackage{overpic}
\usepackage{multirow}
\usepackage{longtable}
\usepackage[dvips]{hyperref}
\usepackage[hang]{subfigure}
\usepackage{amsmath,amssymb, amsthm}
\usepackage[dvips]{color}
\usepackage{framed}
\usepackage[mathscr]{euscript}
\usepackage{enumerate}
\usepackage[toc,page,title,titletoc,header]{appendix} 
\usepackage{cite}

\usepackage{CJKutf8}

 \newtheorem{algorithm}{Algorithm}

\newcommand{\bbf}{\boldsymbol{f}}
\newcommand\bI{\boldsymbol{I}}
\newcommand\bJ{\boldsymbol{J}}
\newcommand\bR{\boldsymbol{R}}
\newcommand\bbR{\mathbb{R}}
\newcommand\bbN{\mathbb{N}}
\newcommand\bxi{\boldsymbol{\xi}}
\newcommand\bx{\boldsymbol{x}}

\newcommand\bq{\boldsymbol{q}}
\newcommand\bu{\boldsymbol{u}}
\newcommand\bv{\boldsymbol{v}}

\newcommand\bLambda{{\boldsymbol{\Lambda}}}

\newcommand\bF{\boldsymbol{F}}

\newcommand\dd{\,\mathrm{d}}
\newcommand\He{\mathit{He}}
\newcommand\Kn{\mathit{Kn}}

\newcommand\mH{\mathcal{H}}

\newcommand\NRxx{NR$xx$~}

\newcommand\mF{\mathcal{F}}

\newcommand\sss{\scriptscriptstyle}
\newcommand\comment[1]{}

\graphicspath{{images/}}
{\theoremstyle{remark} }

\title{A Nonlinear Multigrid Steady-State Solver for Microflow}

\author{ Zhicheng Hu\thanks{LMAM \& School of Mathematical Sciences,
    Peking University, Beijing, China, email: {\tt
      huzhicheng1986@gmail.com}.} ~and Ruo Li\thanks{HEDPS \& CAPT,
    LMAM \& School of Mathematical Sciences, Peking University,
    Beijing, China, email: {\tt rli@math.pku.edu.cn}.}  }

\begin{document}
\maketitle

\begin{abstract}
  We develop a nonlinear multigrid method to solve the steady state of
  microflow, which is modeled by the high order moment system derived
  recently for the steady-state Boltzmann equation with ES-BGK
  collision term. The solver adopts a symmetric Gauss-Seidel iterative
  scheme nested by a local Newton iteration on grid cell level as its
  smoother. Numerical examples show that the solver is insensitive to
  the parameters in the implementation thus is quite robust. It is
  demonstrated that expected efficiency improvement is achieved by the
  proposed method in comparison with the direct time-stepping scheme.

\vspace*{4mm}

\noindent {\bf Keywords:} Multigrid; Boltzmann equation; ES-BGK model;
\NRxx method 
\end{abstract}

\section{Introduction} 
\label{sec:intro} 
The Boltzmann equation plays an important role in various fields of
modern kinetic theory, e.g., rarefied gas flows, microflows and
semiconductor device simulations. Due to its high dimension of
variables, numerical simulation of these problems is extremely
expensive, even though its complicated collision operator (see
e.g. \cite{Cowling}) is replaced by other simplified collision
operators, such as the Bhatnagar-Gross-Krook (BGK) model \cite{BGK},
the ellipsoidal statistical BGK (ES-BGK) model \cite{Holway}, the
Shakhov model \cite{Shakhov}, the Maxwell molecules model
\cite{Ernst}, etc. Lots of work has been done to reduce the
computational costs in numerical solution of the Boltzmann equation or
its simplified collision models over the past few decades. The moment
method, originally introduced by Grad \cite{Grad}, was considered as
one of the candidates in this direction.

The objective of the moment method is to approximate the Boltzmann
equation using a small number of variables. Based on the Hermite
expansion in the velocity space of the distribution function, moment
equations were derived from the Boltzmann equation. These equations
can be viewed as extended hydrodynamic-like equations in a macroscopic
point of view. In terms of numerical method, they are actually a
semi-discretization of the Boltzmann equation, where the velocity
space is discretized by the Hermite spectral method. Moreover, these
equations are highly nonlinear and coupled with each other, as the
classical Euler equations. Meanwhile, the lack of global hyperbolicity
of Grad's original moment system \cite{Muller, Grad13toR13} restricts
its application during a long time. Recently, with a certain
regularization, a systematic numerical method, abbreviated as the
\NRxx method, has been developed in \cite{NRxx,NRxx_new,Cai,Li} to
solve the moment system derived for the Boltzmann equation. The
unified framework of the \NRxx method makes it is easy to implement
the system with moments up to arbitrary order. In \cite{Fan} and
\cite{Fan_new}, a new regularization without any additionally
empirical parameters was proposed such that the resulting moment
system is globally hyperbolic, which is essential in the local
well-posedness of the system.

While the \NRxx method has been found to be effective for various
problems, see e.g. \cite{Wang, Ruo2012, Microflows1D, Hu2012}, the
state-of-the-art numerical techniques to improve its efficiency have
not been sufficiently explored, especially for the steady-state
problems. On the other hand, there are quite some important
applications for the Boltzmann equation, where the main concern is the
steady-state solution. For steady-state problems, some additional
acceleration techniques, including implicit time-stepping schemes
\cite{Mieussens, Mieussens2004}, various iterative methods
\cite{mavriplis1998convergence}, multigrid accelerations
\cite{brandt2011book, hackbusch1985book}, which are originally
developed for classical hydrodynamics \cite{li2008multigrid,
  hu2010robust, hu2011robust, mavriplis2002assessment,
  mavriplis2006multigrid}, may be employed to further improve the
numerical efficiency. Aiming on improving the efficiency of the
steady-state solver following these methods, we are mainly concerned
in this paper with the development of multigrid solution strategy for
the moment system derived for the steady-state Boltzmann equation.

The basic procedure of our current exploration is as below. At first,
we discretize the target moment system under the unified framework of
the \NRxx method as in \cite{NRxx}, so that the multigrid solver we
develop is also unified for the system with moments up to arbitrary
order. We then present a nonlinear iterative method on a single grid
level to solve the resulting discretized system, which is a set of
nonlinear coupled equations as mentioned in the following
sections. This nonlinear iteration is formulated as a nested iterative
scheme which is combined by an inner iteration and an outer
iteration. We take a cell by cell symmetric Gauss-Seidel (SGS)
iteration as the outer iteration to reduce the global nonlinear system
into a local nonlinear system with respect to the local variables on
each cell. For the inner iteration, the Newton type method is employed
to solve the local nonlinear system on each cell. Due to the nonlinear
coupling of the variables in the system, it is quite difficult to
derive the analytical expression of the local Jacobian matrix. As an
alternative, we calculate the local Jacobian matrix through the
numerical differentiation and regularize it by the local residual as
in \cite{li2008multigrid}. The relaxation parameter, in the step of
updating solution of the Newton iteration, is computed adaptively to
preserve positive density and temperature. It is verified by numerical
examples that our method converges much faster than the direct
time-stepping \NRxx scheme.

The nonlinear iteration is taken as smoother in a multigrid framework
for further acceleration. The multigrid method we adopt is the
nonlinear multigrid (NMG) algorithm developed in
\cite{hackbusch1985book}. The NMG algorithm therein is quite formal,
therefore we only need a strategy to generate a sequence of coarse
grids from the finest grid, and operators to transfer the solution
between two successive grids. For the moment system under our current
consideration, which is 1D in spatial coordinates, the coarse grid can
be generated directly by an intuitive way, namely, a coarse grid cell
consists of two adjacent fine grid cells. After the coarse grids have
been obtained, the restriction operator, to transfer the fine grid
solution onto the coarse grid, is then constructed locally. Precisely,
the coarse grid solution on a coarse grid cell is only determined by
the fine grid solution on the corresponding fine grid cells. For the
prolongation operator from the coarse grid to the fine grid, the
identical operator is simply utilized.

The rest part of this paper is organized as follows. In section
\ref{sec:model}, we give a brief review of the steady-state Boltzmann
equation with ES-BGK collision term and its hyperbolic moment system,
followed with the unified discretization. We then present in section
\ref{sec:smoother} the details of the basic nonlinear iterative method
on a single grid level. In section \ref{sec:multi}, the nonlinear
multigrid solver using the basic iteration as smoother is
introduced. Two numerical examples are carried out in section
\ref{sec:example} to demonstrate the robustness and efficiency of the
proposed multigrid solver. Some concluding remarks are given in the
last section.


\section{Moment System of Boltzmann Equation}
\label{sec:model}
The Boltzmann equation in steady state can be written as
\begin{equation}
  \label{eq:boltzmann}
  \bxi \cdot \nabla_{\bx} f + \bF \cdot \nabla_{\bxi} f = Q(f),
\end{equation}
where $f(\bx,\bxi)$ is the particle distribution function of position
$\bx\in \Omega \subset \bbR^D$ $(D=1,2, \text{or } 3)$ and molecular
velocity $\bxi\in \bbR^3$. The vector $\bF$ is the external force
accelerating particles, and the right hand side $Q(f)$ is the
collision term modeling the interaction between particles. Since the
original Boltzmann collision term contains a five-dimensional integral
(see e.g. \cite{Cowling} for the detailed form), it turns out to be
too complicated to handle for numerical solution. Therefore, a variety
of simplified collision models were raised to approximate it.  In this
paper, we take the ES-BGK model as an example to present our method.

The ES-BGK collision term reads
\begin{equation}
  \label{eq:collision-ES}
  Q(f) =  \nu ( f^{\text{ES}} - f),
\end{equation}
where $\nu$ is the average collision frequency given by
\begin{equation}
  \label{eq:nu}
  \nu = \Pr \frac{\rho \theta}{\mu}, 
\end{equation}
and $f^{\text{ES}}$ is an anisotropic Gaussian distribution with the
form as
\begin{equation}
  \label{eq:ES-dis}
  f^{\text{ES}}(\bx,\bxi) = \frac{\rho(\bx)}{m \sqrt{ \det[2\pi \bLambda(\bx)]}
  } \exp\left(-\frac{1}{2} ( \bxi - \bu(\bx) )^T [\bLambda(\bx)]^{-1}
    (\bxi-\bu(\bx)) \right),
\end{equation}
in which $\bLambda=(\lambda_{ij})$ is a $3\times 3$ matrix with
\begin{align}
  \lambda_{ij}(\bx) = \theta(\bx) \delta_{ij} + \left( 1-
    \frac{1}{\Pr} \right) \frac{\sigma_{ij}(\bx)}{\rho(\bx)}, \quad i,
  j = 1,2,3.
\end{align}
Here, $\Pr$ is the Prandtl number, $\mu$ is the viscosity, and $m$ is
the mass of a single particle. The macroscopic quantities, including
density $\rho$, velocity $\bu$, temperature $\theta$, and stress
tensor $\sigma_{ij}$, are related with $f$ by
\begin{align}
  \label{eq:moments}
\begin{aligned}
  & \rho(\bx) = m \int_{\bbR^3} f(\bx,\bxi) \dd \bxi, \quad \rho(\bx)
  \bu(\bx) = m \int_{\bbR^3} \bxi f(\bx,\bxi) \dd \bxi, \\ & \rho(\bx)
  \vert \bu(\bx) \vert^2 + 3 \rho(\bx) \theta(\bx) = m \int_{\bbR^3}
  \vert \bxi \vert^2 f(\bx, \bxi) \dd \bxi, \\ & \sigma_{ij}(\bx) = m
  \int_{\bbR^3} (\xi_i - u_i(\bx))(\xi_j - u_j(\bx)) f(\bx,\bxi) \dd
  \bxi - \rho(\bx) \theta(\bx) \delta_{ij}, \quad i,j = 1,2,3.
\end{aligned}
\end{align}
Additionally, the heat flux $\bq$ is defined by
\begin{equation}
  \label{eq:heat_flux}
  \bq(\bx) = \frac{m}{2} \int_{\bbR^3} \vert \bxi-\bu(\bx) \vert^2
  (\bxi-\bu(\bx)) f(\bx,\bxi) \dd \bxi.
\end{equation}

In \cite{Microflows1D}, a hyperbolic moment system was derived for the
time-dependent Boltzmann equation with the ES-BGK collision term and
without the external force term (see \cite{Li} for the treatment of
the force term). By setting the time derivatives to 0 in that system,
we can obtain a steady-state counterpart for equation
\eqref{eq:boltzmann}. Nevertheless, we briefly review the derivation
of the steady-state moment system below.

The starting point is the approximation of $f$ by an $M$-th order
truncated series as
\begin{equation}
  \label{eq:trucated-dis}
  f(\bx, \bxi) \approx \sum_{\vert \alpha \vert \leq M} f_\alpha(\bx)
  \mH_{\theta(\bx), \alpha} \left( \frac{\bxi - \bu(\bx)}{\sqrt{\theta(\bx)}} \right),
\end{equation}
where $M>2$ is a positive integer, $\alpha =
(\alpha_1,\alpha_2,\alpha_3) \in \bbN^3$ is a three-dimensional
multi-index and $|\alpha| = \alpha_1 + \alpha_2 + \alpha_3$. The
basis functions $\mH_{\theta,\alpha}$ are defined by
\begin{equation}
  \label{eq:base}
  \mH_{\theta,\alpha}(\bv) = \frac{1}{m(2\pi \theta)^{3/2}
    \theta^{|\alpha|/2}} \prod\limits_{d=1}^3 \He_{\alpha_d}(v_d)\exp
  \left(-\frac{v_d^2}{2} \right), \quad \forall \alpha \in
  \bbN^3,~\bv\in\bbR^3,
\end{equation}
where $\He_{n}(x)$ is the Hermite polynomial of order $n$ as
\begin{equation}
  \label{eq:hermite}
  \He_n(x) = (-1)^n\exp \left( \frac{x^2}{2} \right) \frac{\dd^n}{\dd x^n} 
  \exp \left(-\frac{x^2}{2} \right).
\end{equation}
With the help of the orthogonality of Hermite polynomials, the
expansion \eqref{eq:trucated-dis} together with \eqref{eq:moments} and
\eqref{eq:heat_flux} yields
\begin{align}\label{eq:moments-relation}
  \begin{aligned}
    & f_0 = \rho, \qquad f_{e_1} = f_{e_2} = f_{e_3} = 0, \qquad
    \sum_{d=1}^3 f_{2e_d} = 0,
    \\ 
    & \sigma_{ij} = (1+\delta_{ij}) f_{e_i+e_j}, \quad q_i = 2
    f_{3e_i} + \sum_{d=1}^3 f_{2e_d+e_i}, \qquad i,j=1,2,3,
\end{aligned}
\end{align}
where $\delta_{ij}$ is the Kronecker delta symbol, and
$e_1$, $e_2$, $e_3$ are introduced to denote the multi-indices
$(1,0,0)$, $(0,1,0)$, $(0,0,1)$, respectively.

Substituting the expansion \eqref{eq:trucated-dis} into the Boltzmann
equation \eqref{eq:boltzmann}, matching the coefficients of the same
basis function, and regularizing it by the regularization proposed in
\cite{Fan_new}, we obtain the final moment system as follows
\begin{equation}
  \label{eq:mnt-eq}
  \begin{split}
    & \sum_{j=1}^D \Bigg[ \left( \theta \frac{\partial f_{\alpha -
          e_j}}{\partial x_j} + u_j \frac{\partial
        f_{\alpha}}{\partial x_j} + (1-\delta_{|\alpha|,M})(\alpha_j +
      1) \frac{\partial f_{\alpha+e_j}}{\partial x_j} \right) \\ &+
    \sum_{d=1}^3 \frac{\partial u_d}{\partial x_j} \left( \theta
      f_{\alpha-e_d-e_j} + u_j f_{\alpha-e_d} +
      (1-\delta_{|\alpha|,M}) (\alpha_j + 1) f_{\alpha-e_d+e_j}
    \right) \\ &+ \frac{1}{2} \frac{\partial \theta}{\partial x_j}
    \sum_{d=1}^3 \left( \theta f_{\alpha-2e_d-e_j} + u_j
      f_{\alpha-2e_d} + (1-\delta_{|\alpha|,M}) (\alpha_j + 1)
      f_{\alpha-2e_d+e_j} \right) \Bigg] \\ &= \sum_{d=1}^3 F_d
    f_{\alpha-e_d} + \nu (f^{\text{ES}}_\alpha - f_\alpha), \qquad
    |\alpha| \leq M,
    \end{split}
\end{equation}
where $f^{\text{ES}}_\alpha$ are coefficients of the expansion of the
Gaussian distribution $f^{\text{ES}}$, namely,
\begin{align}
  \label{eq:dis-ES-expansion}
  f^{\text{ES}}(\bx, \bxi) = \sum_{|\alpha| \leq M}
  f^{\text{ES}}_{\alpha}(\bx) \mH_{\theta(\bx),\alpha} \left(
    \frac{\bxi - \bu(\bx)}{\sqrt{ \theta(\bx) }} \right).
\end{align}
It follows from \cite{Microflows1D} that $f^{\text{ES}}_\alpha$ can be
calculated recursively by
\begin{align}
  \label{eq:ES-moments}
  f^{\text{ES}}_{\alpha} = \left\{ 
    \begin{aligned}
      & \rho, & & \text{if}~ \alpha = 0, \\
      & 0, & & \text{if}~ |\alpha| = 1, \\
      & \frac{1-1/\Pr}{\alpha_i \rho} \sum_{j=1}^3 \sigma_{ij}
      f^{\text{ES}}_{\alpha - e_i - e_j}, & & {\text{if}}~ |\alpha|
      \geq 2 ~{\text{and}}~ \alpha_i > 0.
    \end{aligned} 
  \right.
\end{align}
At first glance, it is sufficient to note that the moment system
\eqref{eq:mnt-eq} is a set of nonlinear coupled equations for the
moments $\bu, \theta, f_\alpha$.

For practical applications, the moment system \eqref{eq:mnt-eq} has to
be equipped with proper boundary conditions. In kinetic theory, the
Maxwell boundary condition proposed in \cite{Maxwell} is frequently
used. The version of the Maxwell boundary condition for moment system
was derived in \cite{Li} and demonstrated its well-posedness for the
hyperbolic moment system in \cite{Microflows1D}. Accordingly, we adopt
these boundary conditions in this paper. We refer to \cite{Li,
  Microflows1D} for details on these boundary conditions.

Since we are focusing on the iterative method to the steady-state
problem, we discretize the steady-state moment system following the
method in \cite{NRxx, Microflows1D}. The distribution function
\eqref{eq:trucated-dis} is the unknown in the discretized problem,
which is constructed by $\bu, \theta$ and $f_\alpha$. For simplicity of
notations, we introduce $\mF_M(\bu,\theta)$ to denote the linear space
spanned by $\mH_{\theta,\alpha}\left( \frac{\bxi - \bu}{\sqrt{\theta}}
\right)$, $|\alpha|\leq M$. It follows that $\mF_M(\bu,\theta)$ forms
a finite dimensional subspace of $L^2\left(\bbR^3,
\exp\left(|\bxi-\bu|^2/(2\theta)\right)\dd \bxi\right)$. In our
numerical scheme, it is also allowed to approximate the distribution
function $f$ in another linear space $\mF_M(\bu',\theta')$ with the
relation \eqref{eq:moments-relation} violated. However, if it is not
pointed out, $f$ always belongs to the space where
\eqref{eq:moments-relation} holds, i.e., the parameters of the space,
$\bu$ and $\theta$, are macroscopic velocity and temperature of $f$
respectively. A fast transformation between two spaces,
$\mF_M(\bu,\theta)$ and $\mF_M(\bu',\theta')$, can be found in
\cite{NRxx}.

We restrict ourselves to one spatial dimensional case in the
following. Suppose the spatial domain, which is an interval $[a, b]$,
is divided by the mesh
\begin{align*}
  a = x_0 < x_1 < \cdots < x_{N-1} < x_N = b,
\end{align*}
and let $\Delta x_i = x_{i+1}-x_i$. Then the finite volume
discretization for the steady-state moment system \eqref{eq:mnt-eq}
over the $i$-th cell $[x_i, x_{i+1}]$ is given by
\begin{align}
  \label{eq:mnt-eq-dis}
  \frac{F(f_i(\bxi),f_{i+1}(\bxi)) - F(f_{i-1}(\bxi),
    f_i(\bxi))}{\Delta x_i} = G(f_i(\bxi)) + Q(f_i(\bxi)),
\end{align}
where $f_i(\bxi) \in \mF_M(\bu_i,\theta_i)$ is the approximation of
the distribution function on the $i$-th cell. Note that the
distribution function on the ghost cells, $f_{-1}(\bxi)$ and
$f_{N}(\bxi)$, are dependent on $f_{0}(\bxi)$ and $f_{N-1}(\bxi)$,
respectively, for the Maxwell boundary condition.

The numerical flux $F(f_i, f_{i+1})$ is defined on $x_{i+1}$, the
boundary between the $i$-th and $(i+1)$-th cells. To compare the
solution with \cite{Microflows1D}, the same numerical flux, a
non-conservative version of the HLL flux, is considered. We omit its
expression here, since it is enough to know that the numerical fluxes,
$F(f_{i-1}, f_{i})$ and $F(f_{i},f_{i+1})$ of \eqref{eq:mnt-eq-dis},
are approximated in $\mF_M(\bu_i,\theta_i)$, i.e.,
\begin{align}
  \label{eq:flux-expansion}
  \begin{aligned}
    & F(f_{i-1},f_i) = \sum_{\vert \alpha \vert \leq M}
    F_{\alpha}(f_{i-1},f_i) \mH_{\theta_i, \alpha} \left( \frac{\bxi -
        \bu_i}{\sqrt{\theta_i}} \right), \\ & F(f_i,f_{i+1}) =
    \sum_{\vert \alpha \vert \leq M} F_{\alpha}(f_{i},f_{i+1})
    \mH_{\theta_i, \alpha} \left( \frac{\bxi - \bu_i}{\sqrt{\theta_i}}
    \right),
  \end{aligned}
\end{align}
and the computation of these numerical fluxes requires the fast
transformation between $\mF_M(\bu_i,\theta_i)$ and $\mF_M(\bu_{i\pm
  1},\theta_{i\pm 1})$.

Similarly to the numerical flux, the external force term $G(f_i)$ and
the collision term $Q(f_i)$, are also approximated in
$\mF_M(\bu_i,\theta_i)$, i.e.,
\begin{align}
  \label{eq:force-expansion}
  & G(f_i(\bxi)) = \sum_{\vert \alpha \vert \leq M} G_{i,\alpha}
  \mH_{\theta_i, \alpha} \left( \frac{\bxi - \bu_i}{\sqrt{\theta_i}}
  \right),\\
  \label{eq:collision-expansion}
  & Q(f_i(\bxi)) = \sum_{\vert \alpha \vert \leq M}
  Q_{i,\alpha} \mH_{\theta_i, \alpha} \left( \frac{\bxi -
      \bu_i}{\sqrt{\theta_i}} \right),
\end{align}
where $G_{i,\alpha} = \sum_{d=1}^3 F_d f_{i,\alpha-e_d}$ and
$Q_{i,\alpha} = \nu (f_{i,\alpha}^{\text{ES}} - f_{i,\alpha})$.


\section{Single Grid Solver}
\label{sec:smoother}
Define the local residual on the $i$-th cell as
\begin{align}
  \label{eq:residual}
  R_i(f_{i-1}, f_{i}, f_{i+1}) = \frac{F(f_i(\bxi),f_{i+1}(\bxi)) -
    F(f_{i-1}(\bxi), f_i(\bxi))}{\Delta x_i} - G(f_i(\bxi)) -
  Q(f_i(\bxi)).
\end{align}
Then the discretization \eqref{eq:mnt-eq-dis} is re-written as
\begin{align}
  \label{eq:residual-eq-0}
  R_i(f_{i-1}, f_{i}, f_{i+1}) = r_i(\bxi),
\end{align}
where $r_i(\bxi)$ is independent of $f$ and is introduced to give a
slightly more general problem. For \eqref{eq:mnt-eq-dis}, we have
$r_i(\bxi)\equiv 0$. It is clear that the problem
\eqref{eq:residual-eq-0} is nonlinear and we prefer a nested iterative
strategy, which uses a cell by cell SGS iteration as the outer
iteration, and a local Newton iteration as the inner iteration.

Given an approximate solution $f_i^n(\bxi)$, $i=0,1,\ldots, N-1$, the
new approximate solution $f_i^{n+1}(\bxi)$ of an SGS iterative step is
formulated into two loops as
\begin{enumerate}
\item Loop $i$ increasingly from 0 to $N-1$, obtain
  $f_i^{n+\frac{1}{2}}(\bxi)$ by solving
  \begin{align}
    \label{eq:gs-interval}
    R_i(f_{i-1}^{n+\frac{1}{2}},f_i^{n+\frac{1}{2}},f_{i+1}^n)= r_i(\bxi). 
  \end{align}
\item Loop $i$ decreasingly from $N-1$ to $0$, obtain
  $f_i^{n+1}(\bxi)$ by solving
  \begin{align}
    \label{eq:gs-reverse-interval}
    R_i(f_{i-1}^{n+\frac{1}{2}},f_i^{n+1},f_{i+1}^{n+1})= r_i(\bxi). 
  \end{align}
\end{enumerate}
Both \eqref{eq:gs-interval} and \eqref{eq:gs-reverse-interval} are
local nonlinear problems with the distribution function $f_i(\bxi)$ as
the only unknown, and thereby are solved by the Newton
method. Removing the superscripts and the dependence on
$f_{i-1}(\bxi)$, $f_{i+1}(\bxi)$, these problems are abbreviated to
the following form
\begin{align}
  \label{eq:residual-eq}
  R_i(f_{i}) = r_i(\bxi).
\end{align}

Now let $f_i^{(m)}(\bxi)$ is an approximation of $f_i(\bxi)$ and 
re-expand $f_i(\bxi)$ in $\mF_M(\bu_i^{(m)},\theta_i^{(m)})$, that is,
\begin{align}
  \label{eq:trucated-dis-old-base}
  f_i(\bxi) = \sum_{\vert \alpha \vert \leq M} f_{i,\alpha}
  \mH_{\theta_i^{(m)}, \alpha} \left( \frac{\bxi -
      \bu_i^{(m)}}{\sqrt{\theta_i^{(m)}}} \right).
\end{align}
It is trivial to see that $f_i(\bxi)$ is determined by the
coefficients $f_{i,\alpha}$ of \eqref{eq:trucated-dis-old-base}, which
follows that the local residual $R_i(f_i)$ is also determined by
$f_{i,\alpha}$. Consequently, the linearization of \eqref{eq:residual-eq}
by the Newton method is given as
\begin{align}
  \label{eq:linear-pro}
  \begin{aligned}
    \sum_{|k|\leq M} \frac{\delta R_i(f_{i}^{(m)})}{\delta f_{i,k}}
    \Delta f_{i,k}^{(m)} = \tilde{R}_i^{(m)} = r_i(\bxi) -
    R_i(f_i^{(m)}),
  \end{aligned}
\end{align}
where the formal derivatives $\delta R_i / \delta f_{i,k}$ are
calculated by the numerical differentiation method as
\begin{align}
  \frac{\delta R_i(f_{i}^{(m)})}{\delta f_{i,k}} =
  \frac{R_i(\tilde{f}_{i}^{(m),k}) - R_i(f_{i}^{(m)})}{\delta f_{i,k}},
\end{align}
where
\begin{align*}
  \tilde{f}_i^{(m),k}(\bxi) = f_i^{(m)}(\bxi) + \delta
  f_{i,k} \mH_{\theta_i^{(m)}, k}\left( \frac{\bxi -
      \bu_i^{(m)}}{\sqrt{ \theta_i^{(m)} }} \right),
\end{align*}
and $\delta f_{i,k}$ is a small perturbation of the coefficients. As
mentioned in previous section, the residual $R_i(f_i)$, as well as
$r_i(\bxi)$ and its formal derivatives $\delta R_i / \delta f_{i,k}$,
is calculated to a result in $\mF_M(\bu_i^{(m)},\theta_i^{(m)})$. Let
the corresponding coefficients be $R_{i,\alpha}(f_i)$, $r_{i,\alpha}$
and $J_{i,\alpha,k}$, respectively. Then by matching the coefficients
of the same basis function, we can deduce an equivalent linear system
for \eqref{eq:linear-pro} as
\begin{align}
  \label{eq:linear-system}
  \bJ \Delta \bbf_i^{(m)} =  \tilde{\bR}_i^{(m)},
\end{align}
where $\bJ = (J_{i,\alpha,k})$ is the Jacobian matrix, and $\Delta
\bbf_i^{(m)}$, $\tilde{\bR}_i^{(m)}$ are the vectors of $\Delta
f_{i,\alpha}^{(m)}$, $\tilde{R}_{i,\alpha}^{(m)} = r_{i,\alpha} -
R_{i,\alpha}(f_i^{(m)})$, respectively.

Since the linear problem \eqref{eq:linear-system} might be singular,
to make the iteration stable while keeping the convergence speed, it
is regularized by the local residual as in \cite{li2008multigrid},
which yields
\begin{align}
  \label{eq:linear-system-regularization}
  \begin{aligned}
    \left(\lambda \left\Vert \tilde{R}_i^{(m)} \right\Vert \bI + \bJ
    \right) \Delta \bbf_i^{(m)} = \tilde{\bR}_{i}^{(m)},
  \end{aligned}
\end{align}
where $\bI$ is the identity matrix, $\lambda$ is a parameter, and
$\left\Vert \tilde{R}_i^{(m)} \right \Vert$ is a norm of the residual
$\tilde{R}_i^{(m)}$. In our numerical examples, $\lambda$ is
insensitive and is often set as 1. The $L^2$ norm, equipped for the
linear space $\mF_M(\bu_i^{(m)}, \theta_i^{(m)})$, is employed for
$\Vert \cdot \Vert$, which is turned out to be
\begin{align}
  \label{eq:l2-norm-dis}
  \left\Vert \tilde{R}_i^{(m)} \right \Vert = \sqrt{
    \sum_{|\alpha|\leq M} C_\alpha \left| \tilde{R}_{i,\alpha}^{(m)}
    \right|^2},
\end{align}
where $C_\alpha = (2\pi)^{-3/2} \left( \theta_i^{(m)
  }\right)^{-|\alpha|-3} \alpha!$ and $\alpha! = \alpha_1!  \alpha_2!
\alpha_3!$.

While obtaining $\Delta \bbf_i^{(m)}$ from
\eqref{eq:linear-system-regularization}, the solution is then updated
by
\begin{align}
  \label{eq:update-dis}
  f_{i,\alpha}^{*} = f_{i,\alpha}^{(m)} + \tau \Delta
  f_{i,\alpha}^{(m)},
\end{align}
where $\tau$ is a relaxation parameter set as $\tau = \min\{1,
\hat{\tau}\}$, in which $\hat{\tau}>0$ is a parameter to preserve the
positivity of the local density and temperature. The computation of
$\hat{\tau}$ is not difficult but tedious, hence we present it in
Appendix \ref{sec:def-tau}.

The updated approximation $f_i^*(\bxi)$, obtained until now, belongs to
$\mF_M(\bu_i^{(m)}, \theta_i^{(m)})$ and usually does not satisfy the
relation \eqref{eq:moments-relation}. An additional step is utilized
such that the expansion of the new approximation satisfies
\eqref{eq:moments-relation}. To this end, compute the new macroscopic
velocity $\bu_i^{(m+1)}$ and temperature $\theta_i^{(m+1)}$ from
$f_i^{*}(\bxi)$, then project $f_i^{*}(\bxi)$ into
$\mF_M(\bu_i^{(m+1)}, \theta_i^{(m+1)})$ to obtain
$f_i^{(m+1)}(\bxi)$.

Finally, the Newton method for the local problem
\eqref{eq:residual-eq} is collected as follows
\begin{algorithm}[local Newton iteration]~\label{alg:local-newton}
\begin{enumerate}
  \item Given an initial guess $f_i^{(0)}(\bxi)$, let $m:=0$.
  \item If $i=0$ or $i=N-1$, compute the distribution function on the
    ghost cell based on the boundary conditions.
  \item Compute the residual $\tilde{R}_i^{(m)}$ and its $L^2$
    norm. If $\left\Vert \tilde{R}_i^{(m)}\right\Vert \leq
    \text{\it{tol}}$, stop; otherwise, go to the next step.
  \item Compute the Jacobian matrix, and solve
    \eqref{eq:linear-system-regularization} for $\Delta
    f_{i,\alpha}^{(m)}$.
  \item Compute the relaxation parameter $\tau$, and update the
    solution by \eqref{eq:update-dis}.
  \item Compute $\bu_i^{(m+1)}$, $\theta_i^{(m+1)}$, and project the
    new approximation into $\mF_M(\bu_i^{(m+1)}, \theta_i^{(m+1)})$ to
    get $f_i^{(m+1)}(\bxi)$.
  \item $m:=m+1$, return to step 2.
  \end{enumerate}
\end{algorithm}
The parameter $\text{\it{tol}}$ in the above algorithm is a local criterion to
determine whether the local steady state is reached. It is usually set
as $\text{\it{Tol}}$, which is used in Algorithm \ref{alg:outest-loop} as the
criterion of the global steady state.

Now the basic nonlinear iteration for the discretized system
\eqref{eq:residual-eq-0} has been obtained. In the rest of this paper,
it is denoted by $f^{\text{new}} = \text{SGS-Newton}(f^{\text{old}},
r, \nu)$, where $f$ is the approximation of the distribution function
with $f_i \in \mF_M(\bu_i, \theta_i)$, $r$ is the right hand side of
\eqref{eq:residual-eq-0}, and $\nu$ is the steps of SGS iteration.

It is sure that the SGS-Newton iteration could be performed until the
global steady state is achieved as the following algorithm.
\begin{algorithm}~\label{alg:outest-loop}
  \begin{enumerate}
  \item Given an initial solution $f^{0}(\bxi)$ with $f_i^{0}\in
    \mF_M(\bu_i^{0}, \theta_i^{0})$, let $n:=0$.
  \item Perform an \emph{SGS-Newton} iteration, i.e., $f^{n+1} =
    \emph{SGS-Newton}(f^{n}, r, 1)$.
  \item Calculate the global residual $\tilde{R}$ with $\tilde{R}_i =
    r_i - R_i(f_{i-1}^{n+1},f_{i}^{n+1},f_{i+1}^{n+1})$, and its $L^2$
    norm, which is defined as
    \begin{align}
      \label{eq:global-l2-norm}
      \left\Vert \tilde{R} \right \Vert = \sqrt{ \sum_{i=0}^{N-1}
        \left\Vert \tilde{R}_{i} \right\Vert^2 \Delta x_i}.
    \end{align}
  \item If $\left\Vert \tilde{R} \right\Vert \leq \text{\it{Tol}}$,
    stop; otherwise, let $n:=n+1$, return to step 2.
  \end{enumerate}
\end{algorithm}

\begin{figure}[!htb]
  \centering
  \includegraphics[width=0.8\textwidth]{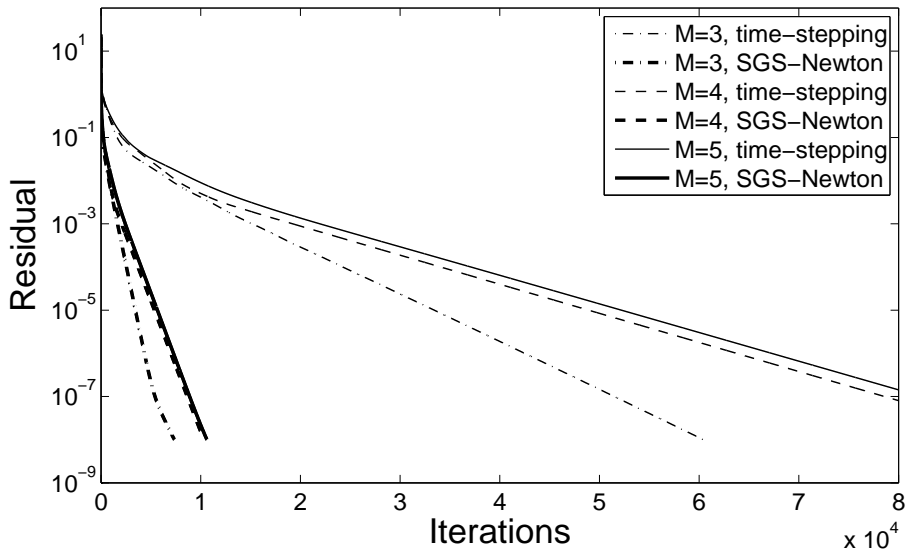}
  \caption{Convergence history on a uniform grid with $N=1024$ of the
    planar Couette flow for $\Kn=0.1199$, $u^W = 1.2577$.}
  \label{fig:comp-sgsn-tsm-res}
\end{figure}

As an efficiency test of the SGS-Newton iteration, we consider the
convergence for the planar Couette flow (see section
\ref{sec:num-ex-couette}) using Algorithm \ref{alg:outest-loop}.
Compared with the explicit time-stepping scheme as in
\cite{Microflows1D}, the convergence history of the SGS-Newton
iteration is shown in Figure \ref{fig:comp-sgsn-tsm-res}. As expected,
the SGS-Newton iteration provides us a faster convergence than the
explicit time-stepping scheme.

In view of that the exact solution for the local problem
\eqref{eq:residual-eq} is not necessary, an improvement in efficiency
is obtained by additional criteria for Algorithm
\ref{alg:local-newton}. In our implementation, the local Newton
iteration also stops if the local residual decreases in half or the
maximum steps, set as 5, is reached. An even essential improvement in
efficiency is multigrid acceleration by using the SGS-Newton iteration
as smoother, which is described in the following section.

\section{Multigrid Solver}
\label{sec:multi}
Using the SGS-Newton iteration as smoother, the multigrid strategy we
adopt is the nonlinear multigrid approach \cite{hackbusch1985book}.
Besides the smoother, the main ingredients of a nonlinear multigrid
solver are the coarse grid correction and the operators (restriction
and prolongation) between the fine and coarse grids. In this section,
we first consider a two-grid solver to give the coarse grid correction
problem, and to construct the restriction and the prolongation
operators. Then it is generalized to a complete multigrid algorithm by
recursion. For convenience, we introduce subscripts $h$ and $H$ to
denote operators and variables related to the fine and coarse grids,
respectively.

\subsection{Coarse grid correction}
\label{sec:multi-coarse}
In contrast to the linear multigrid procedure which solves a
correction equation for the coarse grid correction directly, the
nonlinear multigrid procedure calculate the coarse grid solution
first. As seen in the last section, let us re-written the fine grid
problem resulting from \eqref{eq:residual-eq} into a global form as
\begin{align}
  \label{eq:multi-fine-problem}
  R_h(f_h) = r_h,
\end{align}
and suppose $\bar{f}_h$ is an approximate solution for the above
problem. The corresponding coarse grid problem is given as
\begin{align}
  \label{eq:multi-coarse-problem}
  R_H(f_H) = r_H \triangleq R_H(I_h^H \bar{f}_h) + I_h^H \left(r_h -
    R_h(\bar{f}_h)\right),
\end{align}
where $I_h^H$ is the restriction operator from the fine grid to the
coarse grid, and $R_H(f_H)$ is evaluated analogously to the fine grid
counterpart $R_h(f_h)$. It follows that
\eqref{eq:multi-coarse-problem} can be solved using the SGS-Newton
iteration, and a solution $f_H$ is obtained by employing Algorithm
\ref{alg:outest-loop}. Then the fine grid solution $\bar{f}_h$ is
corrected as
\begin{align}
  \label{eq:multi-update}
  \hat{f}_h = \bar{f}_h + I_H^h\left( f_H - I_h^H \bar{f}_h \right),
\end{align}
where $I_H^h$ is the prolongation operator from the coarse grid to the
fine grid.

\subsection{Restriction and prolongation}
\label{sec:multi-restriction}
In our implementation, the coarse grid is generated from the fine grid
by a standard way, namely, the coarse grid point $x_{\sss H,i}$
coincides with the fine grid point $x_{\sss h, 2i}$. Based on this, we
construct the restriction operator $I_h^H$ locally, which means for
any fine grid function $g_{\sss h}$, its restriction $g_{\sss H,i}$ on
the $i$-th coarse grid cell $[x_{\sss H,i},x_{\sss H,i+1}]$, given by
$g_{\sss H}=I_h^H g_{\sss h}$, is determined only by $g_{\sss h,2i}$
and $g_{\sss h,2i+1}$. The detailed construction of the restriction
operator $I_h^H$ can be found in Appendix \ref{sec:def-restriction},
and we only show the final result here.

There are two grid variables in \eqref{eq:multi-coarse-problem}, the
solution $\bar{f}_h$ and the residual $\bar{R}_h=r_h-R_h(\bar{f}_h)$,
which are required to be transferred on to the coarse grid using
$I_h^H$. The restrictions $\bar{f}_{H,i}$ and $\bar{R}_{H,i}$ on the
$i$-th coarse grid cell are constructed in the following steps
\begin{enumerate}
\item Compute $\bar{\bu}_{\sss H,i}$ and $\bar{\theta}_{H,i}$ from
  \eqref{eq:restriction-rho-u-theta};
\item Project $\bar{f}_{h,2i}$, $\bar{f}_{h,2i+1}$, $\bar{R}_{h,2i}$
  and $\bar{R}_{h,2i+1}$ into
  $\mF_M(\bar{\bu}_{\sss H,i},\bar{\theta}_{H,i})$, then calculate
  $\bar{f}_{H,i,\alpha}$ and $\bar{R}_{H,i,\alpha}$ by the formula
  \eqref{eq:restriction-moments}.
\end{enumerate}

For the prolongation operator $I_H^h$, the simple identical operator
is employed. Then the correction formula \eqref{eq:multi-update} is
re-written as
\begin{align}
  \label{eq:multi-update-simplify}
  \hat{f}_h = \delta \bar{f}_h + f_H,
\end{align}
where $\delta \bar{f}_h = \bar{f}_h - I_h^H \bar{f}_h$. It is pointed
out that the above formula includes the projection procedure to
represent $\hat{f}_h$ in $\mF_M(\hat{\bu}_{\sss h},\hat{\theta}_{h})$.

\subsection{Overall algorithm}
\label{sec:multi-alg}
As the computation of the exact solution $f_{H}$ on the coarse grid by
the SGS-Newton iteration is still cumbersome for a large number of
$N_H$, the coarse grid problem \eqref{eq:multi-coarse-problem} is
solved recursively with the multigrid algorithm until the exact
solution on the coarsest grid can be obtained cheaply. This results in
a complete multigrid algorithm.

Now let us introduce subscripts $h_k$, $k=0,1,\ldots,K$ to denote
operators and variables related to the $k$-th level grid, where $h_0$
and $h_{K}$ correspond to the coarsest and the finest grid,
respectively. Then the $(k+1)$-th level multigrid iteration, denoted by
$f_{h_k}^{n+1} = \text{NMG}(f_{h_k}^n, r_{h_k}, k )$, is given in the
following algorithm.
\begin{algorithm} [Nonlinear multigrid (NMG) iteration]~
\label{alg:nmg}
  \begin{enumerate}
  \item If $k=0$, call Algorithm \ref{alg:outest-loop} to have a
    solution as $f_{h_0}^{n+1}$; otherwise, go to the next step.
  \item Pre-smoothing: perform $\nu_1$ smoothing steps using the
    \emph{SGS-Newton} iteration to obtain a new approximation
    $\bar{f}_{h_k}$, that is, $\bar{f}_{h_k} =
    \emph{SGS-Newton}(f_{h_k}^n, r_{h_k}, \nu_1)$.
  \item Coarse grid correction:
    \begin{enumerate}
    \item Compute the fine grid residual as $\bar{R}_{h_k} = r_{h_k}
      - R_{h_k}(\bar{f}_{h_k})$.
    \item Calculate the coarse grid approximation by the restriction
      operator $I_{h_k}^{h_{k-1}}$ as $\bar{f}_{h_{k-1}} =
      I_{h_k}^{h_{k-1}} \bar{f}_{h_k}$, and compute the difference
      $\delta \bar{f}_{h_k}$ with 
      \begin{align*}
        \delta \bar{f}_{h_k,2i} = \bar{f}_{h_k,2i} -
        \bar{f}_{h_{k-1},i}, \quad \delta \bar{f}_{h_k,2i+1} =
        \bar{f}_{h_k,2i+1} - \bar{f}_{h_{k-1},i}, \qquad
        i=0,1,\ldots,N_{h_{k-1}}.
      \end{align*}
    \item Calculate the right hand side of the coarse grid problem
      \eqref{eq:multi-coarse-problem} as $r_{h_{k-1}} =
      I_{h_k}^{h_{k-1}} \bar{R}_{h_k} + R_{h_{k-1}}(
      \bar{f}_{h_{k-1}})$.
    \item Recursively call the multigrid algorithm (repeat $\gamma$
      times with $\gamma=1$ for a so-called $V$-cycle, $\gamma=2$ for a
      $W$-cycle, and so on) as
      \begin{align*}
        \tilde{f}_{h_{k-1}} = \emph{NMG}^\gamma (\bar{f}_{h_{k-1}},
        r_{h_{k-1}}, k-1).
      \end{align*}
    \item Correct the fine grid solution by $\hat{f}_{h_{k}} = \delta
      \bar{f}_{h_k} + \tilde{f}_{h_{k-1}}$.
    \end{enumerate}
  \item Post-smoothing: perform $\nu_2$ smoothing steps using the
    \emph{SGS-Newton} iteration to obtain the new approximation
    as
    \begin{align*}
      f_{h_k}^{n+1} = \emph{SGS-Newton}(\hat{f}_{h_k}, r_{h_k},
      \nu_2).
    \end{align*}
  \end{enumerate}
\end{algorithm}
The NMG solver is then obtained from Algorithm \ref{alg:outest-loop}
where the SGS-Newton iteration is replaced by the above NMG iteration.


\section{Numerical Examples}
\label{sec:example}
Two numerical examples, the planar Couette flow and the force driven
Poiseuille flow, are presented in this section to illustrate the main
features of the NMG solver described in previous section. In all tests
below, a $V$-cycle NMG with smoothing steps $\nu_1=\nu_2=2$ is
employed, and its level is chosen such that there are only four cells
on the coarsest grid. With such a small grid, the solution can
be obtained very cheaply by the SGS-Newton iteration, e.g., in most
tests, only $30$ SGS-Newton iterations is enough to resolve the
coarsest grid problem exactly at the beginning of the NMG iteration,
and the steps of the SGS-Newton iteration is quickly decreased to 1 in
successive NMG iterations.

For simplicity, the Prandtl number $\Pr$ is set as $2/3$, and we
consider the solution of the dimensionless equation, which has the
same form as \eqref{eq:boltzmann} with the particle mass $m=1$. As
mentioned in Section \ref{sec:model}, the Maxwell boundary conditions
for moment system derived in \cite{Li} are adopted. It is easy to show
that such boundary conditions could not determine a unique solution
for the steady-state moment system \eqref{eq:mnt-eq}. For the
time-stepping scheme, the additional initial condition is used to
obtain the unique solution. However, since the classical Gauss-Seidel
iteration is non-conservative, it is apparently that our NMG iteration
using the SGS-Newton iteration as smoother is also
non-conservative. This leads to a steady-state solution which is
inconsistent with the time-stepping solution only for the density,
even though the same initial condition is employed as the initial
guess in our NMG solver. In order to converge towards the same
steady-state solution as the time-stepping scheme, the solution is
corrected as \cite{Mieussens2004} at each NMG iterative step, that is,
\begin{align}
  \label{eq:preserve-mass}
  f^{\text{new}} \leftarrow \frac{\displaystyle\int f^0(x, \bxi) \dd
    \bxi \dd x}{\displaystyle\int f^{\text{new}}(x,\bxi) \dd \bxi \dd
    x} f^{\text{new}} = \frac{\sum\limits_{i=0}^{N-1} \rho_i^0 \Delta
    x_i}{\sum\limits_{i=0}^{N-1} \rho_i^{\text{new}} \Delta x_i}
  f^{\text{new}}.
\end{align}
The above correction is sufficient to recover the consistent
steady-state solution.

All our computations are performed on Linux operating system on an
Xeon workstation with a quad-core processor and core speed
2.93GHz. The global tolerance $\text{\it{Tol}}$ is set as $10^{-8}$.

\subsection{The planar Couette flow}
\label{sec:num-ex-couette}
We first consider a force-free benchmark problem, namely, the planar
Couette flow. To compare with the solution of the moment system using
the time-stepping scheme in \cite{Microflows1D}, we use the same data
and parameters in our tests. The gas of argon lies between two plates
parallel to the $yz$-plane with a distance of $L=1$. These two plates
have the same temperature of $\theta^W=1$, and one plate is stationary
while the other moves with a constant velocity $u^W$ in the $y$
direction. The dimensionless collision frequency $\nu$ is given as
\begin{align}
  \label{eq:couette-nu}
  \nu = \sqrt{\frac{\pi}{2}} \frac{\Pr}{\Kn} \rho \theta^{1-w},
\end{align}
where $\Kn$ is the Knudsen number, and the viscosity $\mu$ used in
\eqref{eq:nu} is assumed as a function directly proportional to $w$
power of temperature, i.e., $\mu \propto \theta^w$. For the gas of
argon, we have $w=0.81$. These settings are actually the same as in
\cite{Mieussens2004}, whose solution of the non-dimensionless
Boltzmann equation using the discrete velocity method is utilized as a
reference of the solution of the moment system.

\begin{figure}[!tb]
{  \centering
  \subfigure[Density, $\rho$]{\includegraphics[width=0.5\textwidth]{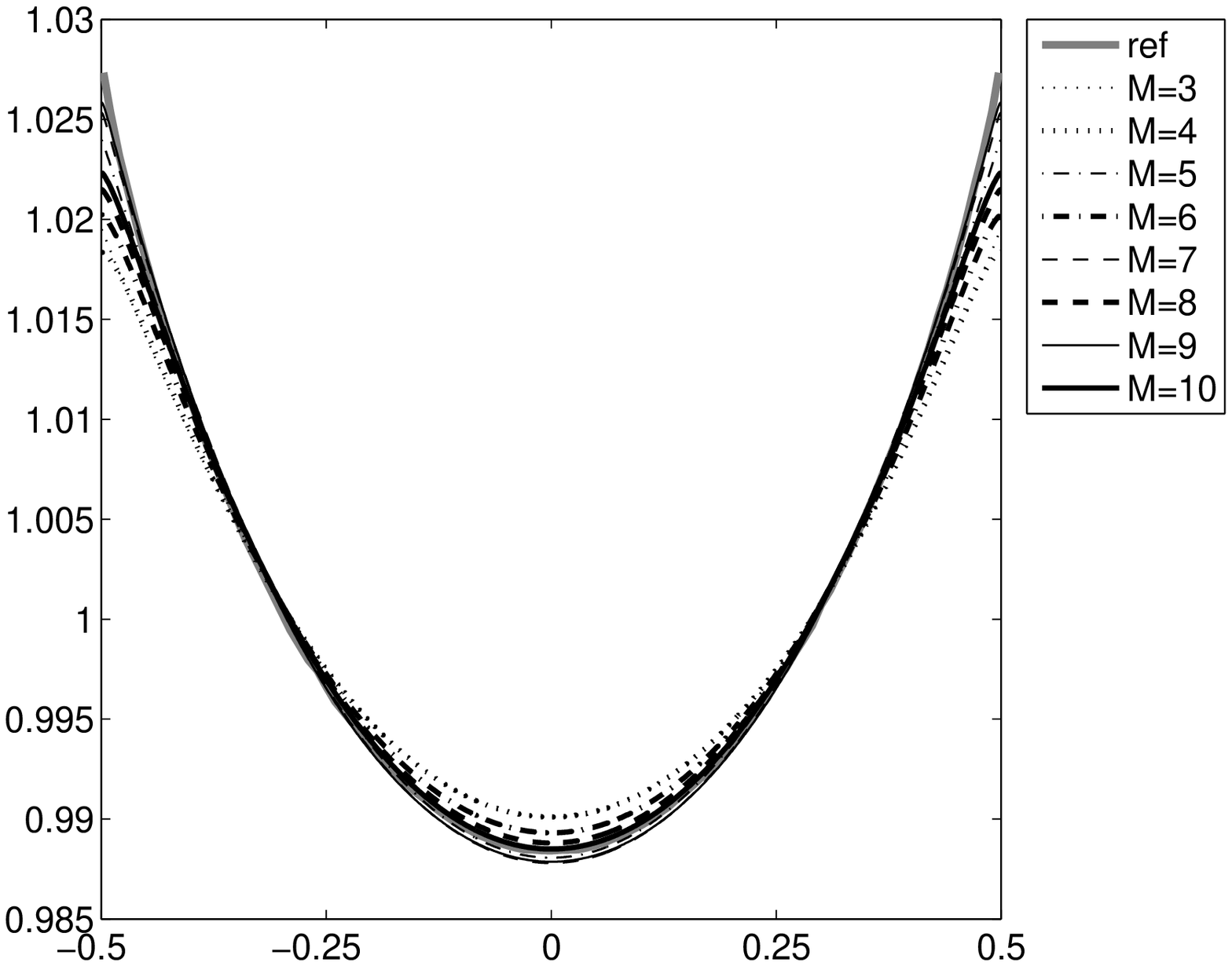}}\hfill
  \subfigure[Temperature, $\theta$]{\includegraphics[width=0.495\textwidth]{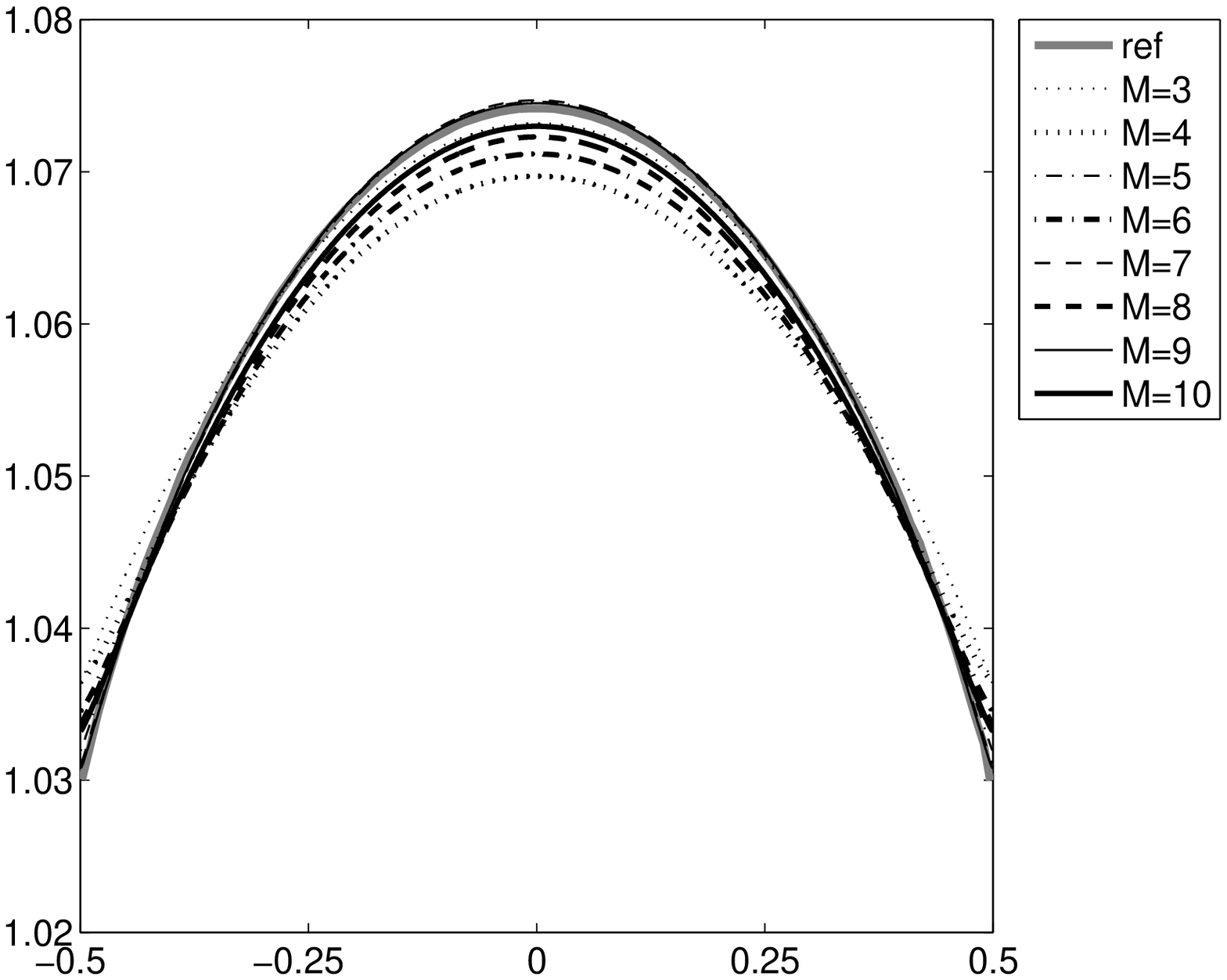}} \\
  \subfigure[Shear stress, $\sigma_{12}$]{\includegraphics[width=0.495\textwidth]{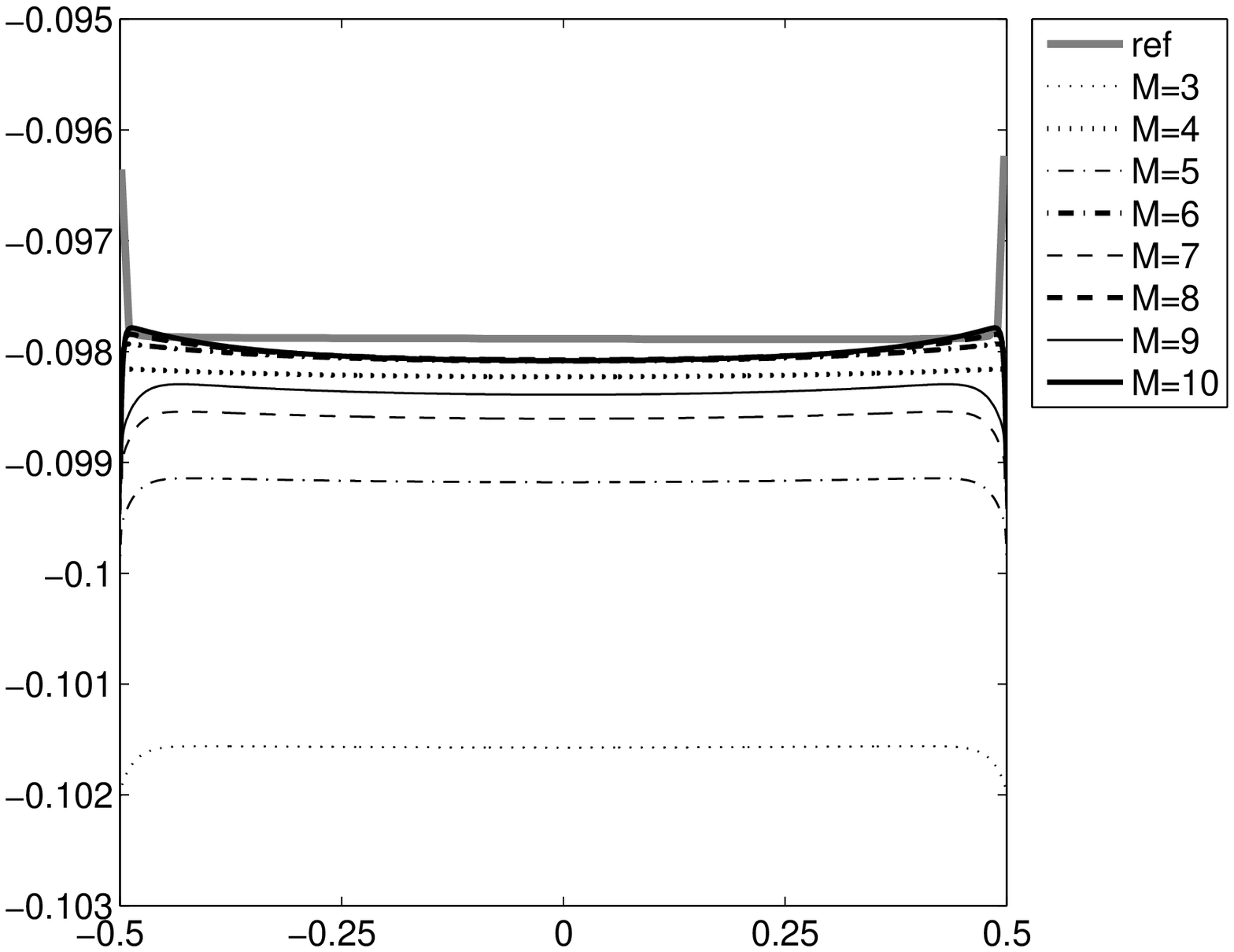}} \hfill
  \subfigure[Heat flux, $q_1$]{\includegraphics[width=0.495\textwidth]{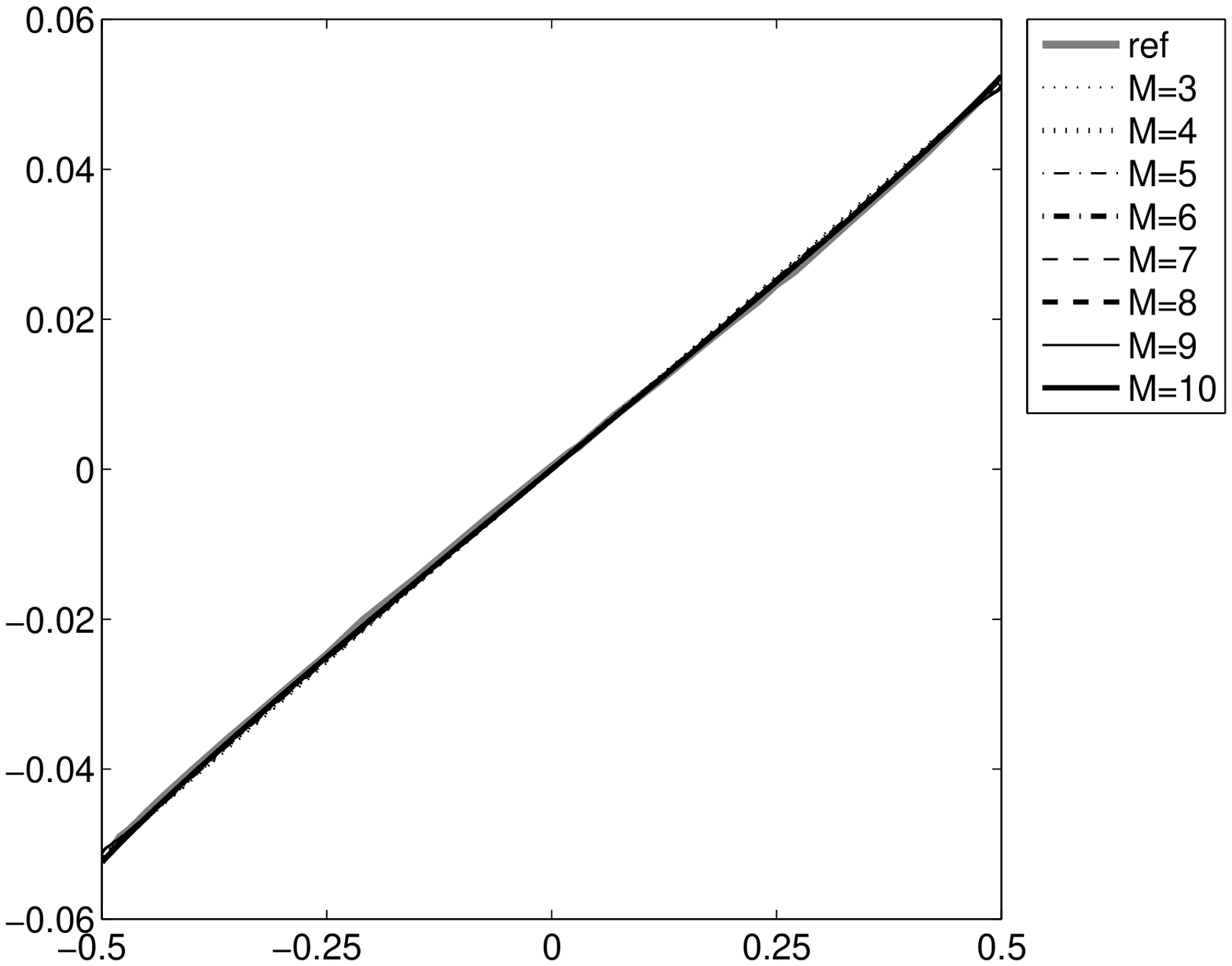}} }
\caption{Solution of the Couette flow for $\Kn=0.1199$ and $u^W =
  1.2577$ on a uniform grid with $N=2048$.}
  \label{fig:couette-Kn01-uw300-density}
\end{figure}

\begin{figure}[!tb]
  \centering
  \includegraphics[width=0.8\textwidth]{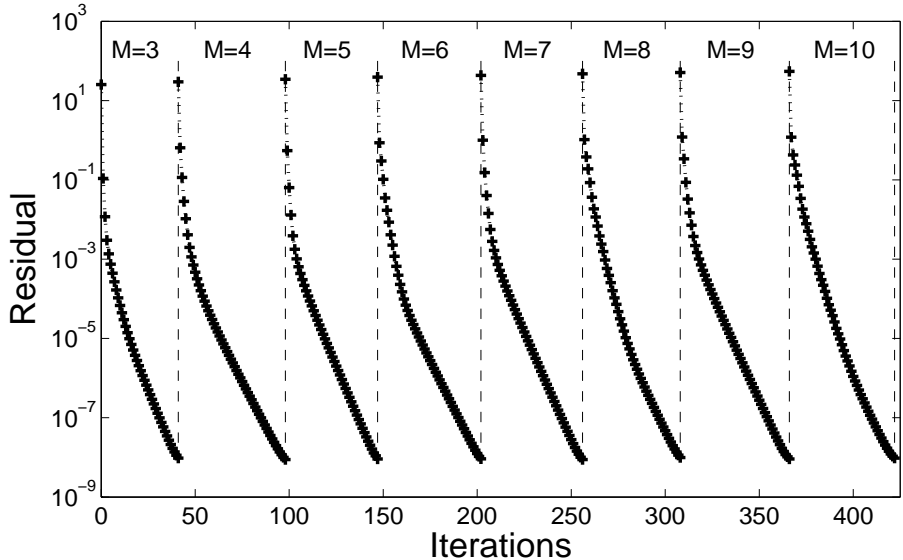}\hfill
  \caption{Convergence history of the Couette flow with $\Kn=0.1199$
    and $u^W = 1.2577$ on a uniform grid with $N=2048$.}
  \label{fig:couette-Kn01-uw300-order-res}
\end{figure}

In this example, all the computations begin with a global equilibrium with
\begin{align}
  \label{eq:couette-initial}
  \rho^0 (x) = 1, \quad \bu^0(x) = 0, \quad \theta^0(x) = 1.
\end{align}
Since our NMG solver delivers exactly the same steady-state solution
as the time-stepping scheme in \cite{Microflows1D}, where the solution
of the moment system for this example has been presented and compared
with the reference solution, we omit any discussion on the accuracy of
our solution, and only focus on the efficiency of the NMG solver. As
an example, the steady-state solution for $\Kn = 0.1199$ and $u^W =
1.2577$ on a uniform grid with $N = 2048$ is displayed in Figure
\ref{fig:couette-Kn01-uw300-density}, compared with the reference
solution. The convergence history in terms of NMG iterations for this
case is presented in Figure \ref{fig:couette-Kn01-uw300-order-res},
which shows the efficiency and robustness of the NMG solver for
various $M$.

\begin{figure}[!tb]
  \centering
  \includegraphics[width=0.5\textwidth]{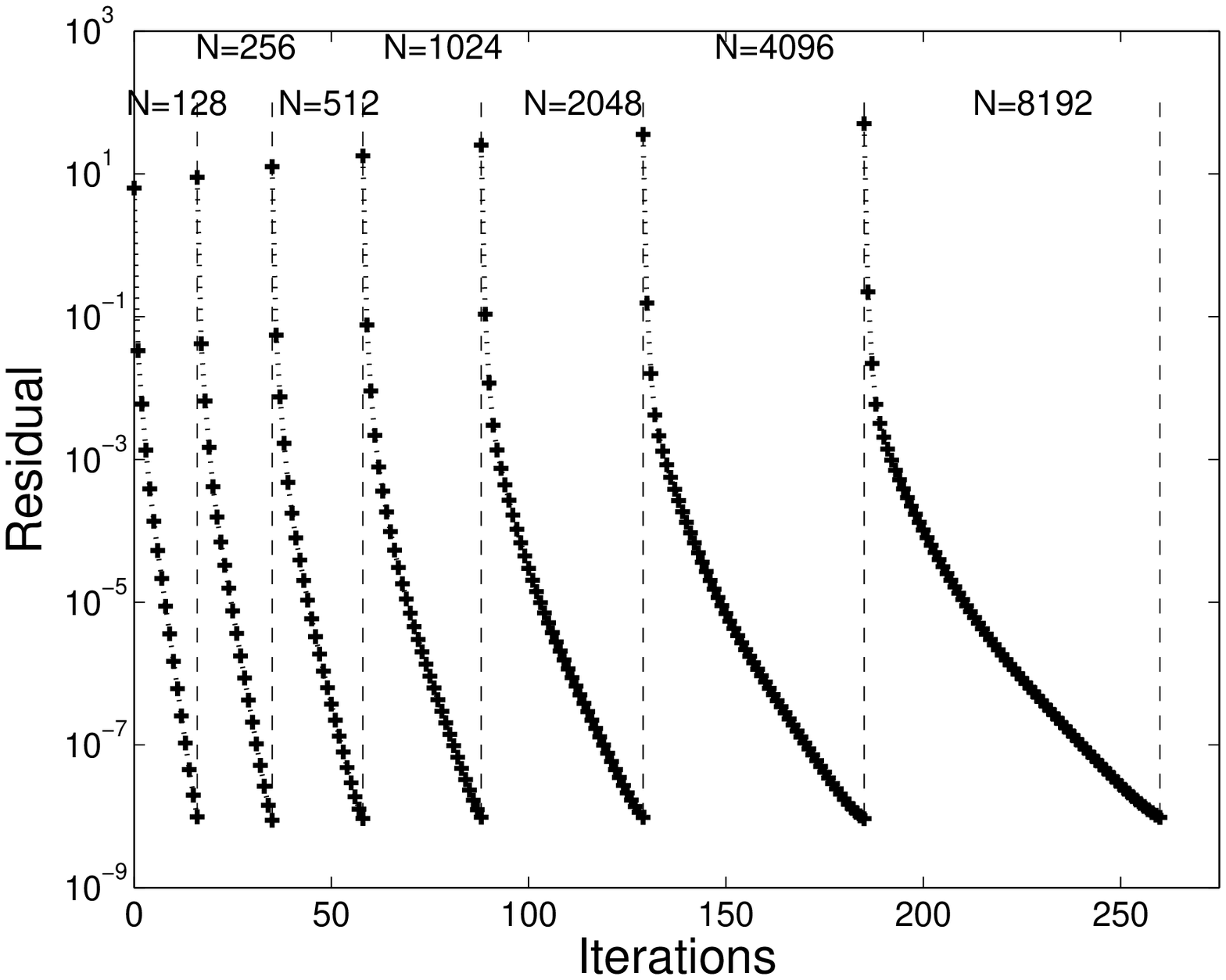}\hfill
  \includegraphics[width=0.5\textwidth]{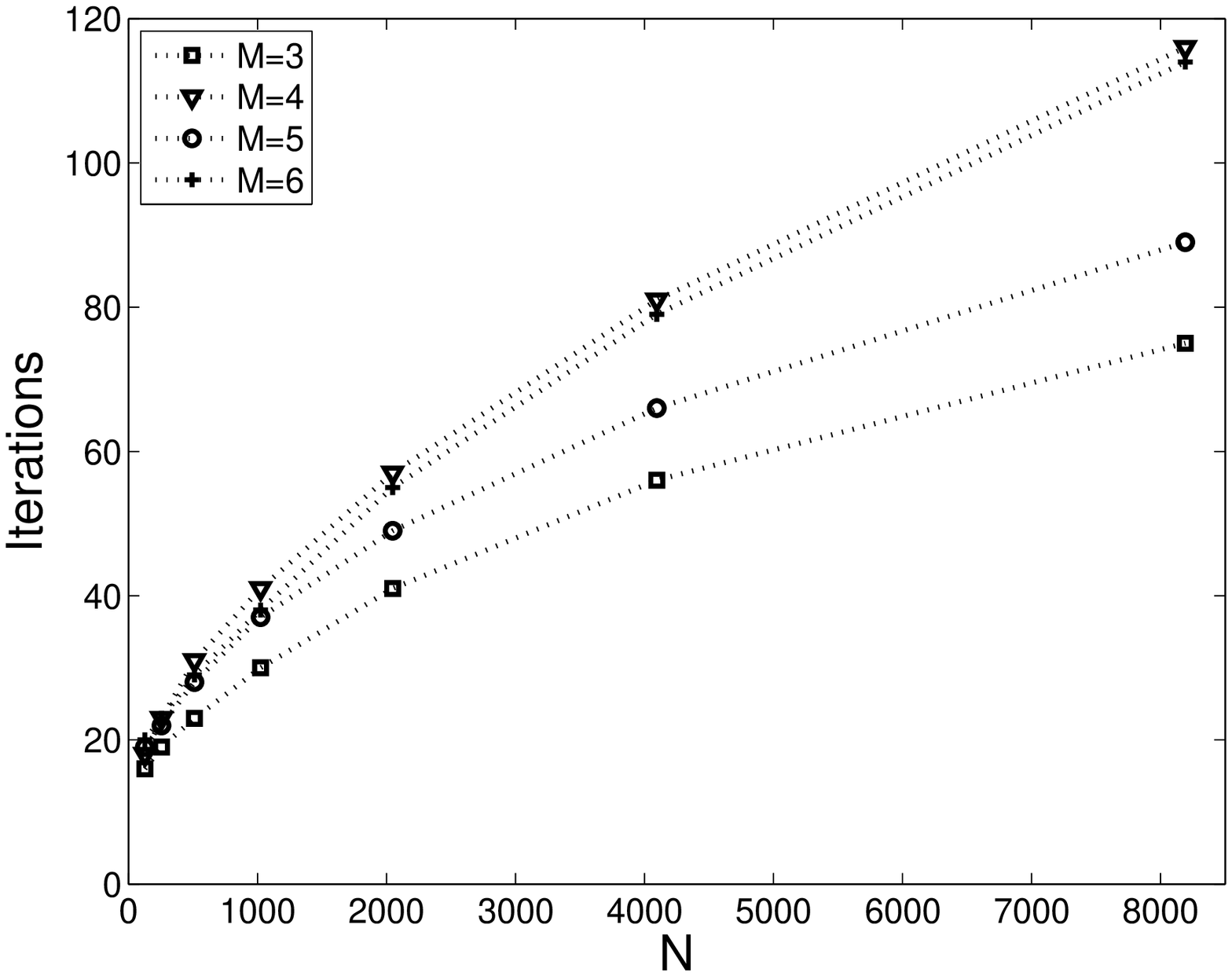}
  \caption{Convergence history ($M=3$, left) and iterations (right) of
    the Couette flow for $\Kn=0.1199$ and $u^W=1.2577$ on different
    uniform grids.}
  \label{fig:couette-Kn01-uw300-mesh-res}
\end{figure}

Now let us consider the behavior of the NMG solver with respect to the
grid size. We first test the Couette flow for $\Kn = 0.1199$ and $u^W
= 1.2577$ on a uniform refined meshes sequence. As shown in Table
\ref{tab:couette-mesh} for the number of NMG iterations, and in Figure
\ref{fig:couette-Kn01-uw300-mesh-res} for convergence history with $M
= 3$ and total iterations in terms of grid size, the NMG solver works
well on the uniform mesh. Although the total iterations increase as
the grid size increases, the rate of increase is much slower than the
single grid solver such as the SGS-Newton iteration and the
time-stepping scheme, which are well known that double the iterations
as the grid size doubles. This indicates that the NMG solver is
substantially more efficient than the single grid solver on finer
grids. It is also shown in Figure
\ref{fig:couette-Kn01-uw300-mesh-res} (right) that the solution of the
moment system for odd $M$ converge faster than the solution for even
$M$ as grid size increases.

\begin{table}[!tb]
  \centering
  {\begin{tabular}{l|l|ccccccc}
 \hline
\multicolumn{2}{c|}{$N$} & $2^7$ & $2^8$ & $2^9$ & $2^{10}$ & $2^{11}$ & $2^{12}$ & $2^{13}$ \\
\hline \multirow{4}{*}{\small Iterations}
& $M=3$ & 16& 19& 23& 30& 41& 56& 75 \\
& $M=4$ & 18& 23& 31& 41& 57& 81& 116 \\
& $M=5$ & 19& 22& 28& 37& 49& 66& 89 \\
& $M=6$ & 20& 23& 29& 38& 55& 79& 114 \\
\hline
  \end{tabular}}
\caption{Iterations of the Couette flow for $\Kn = 0.1199$ and $u^W = 1.2577$ on different uniform grids.}
  \label{tab:couette-mesh}
\end{table}

\begin{figure}[!tb]
\includegraphics[width=0.5\textwidth]{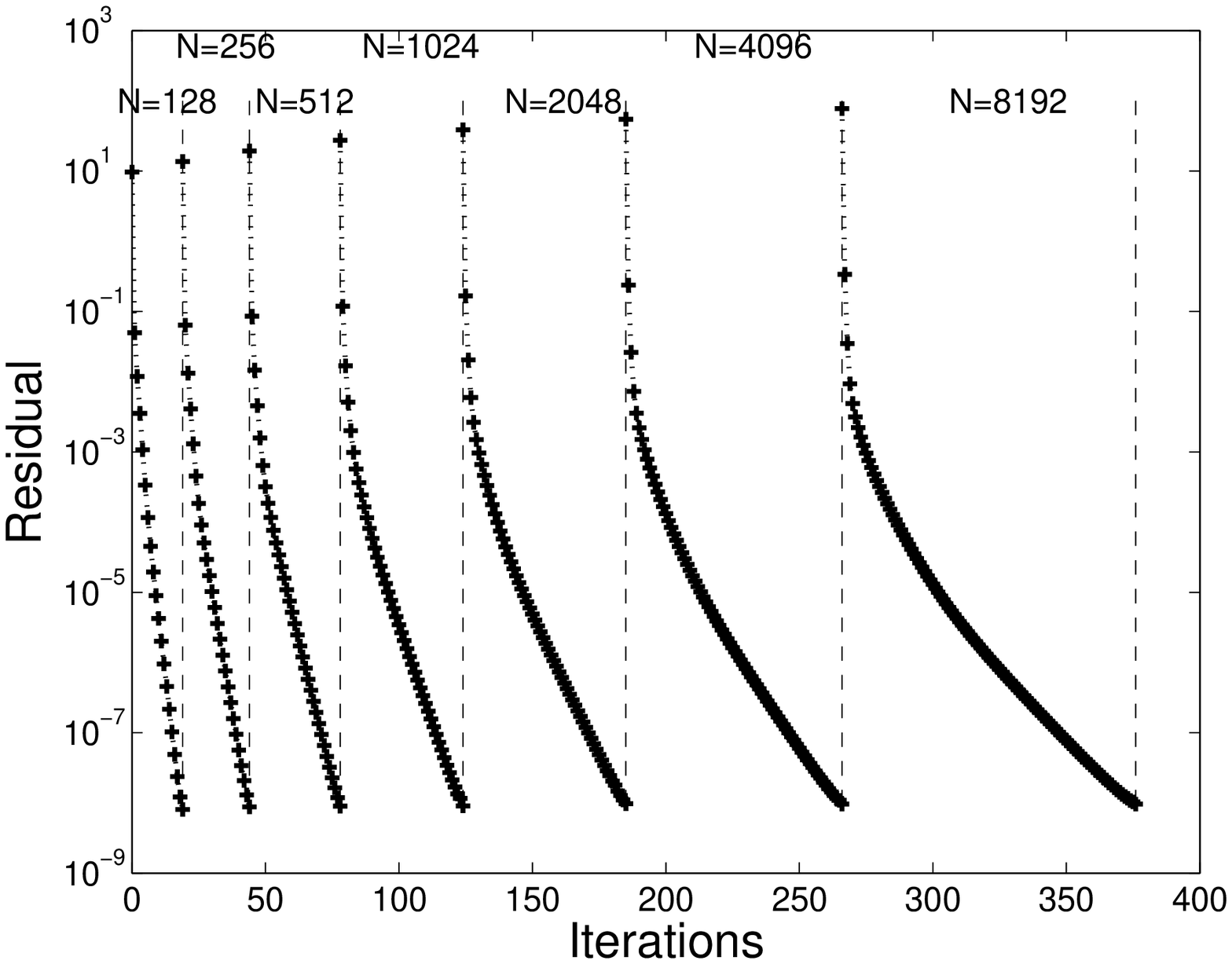}
  \includegraphics[width=0.5\textwidth]{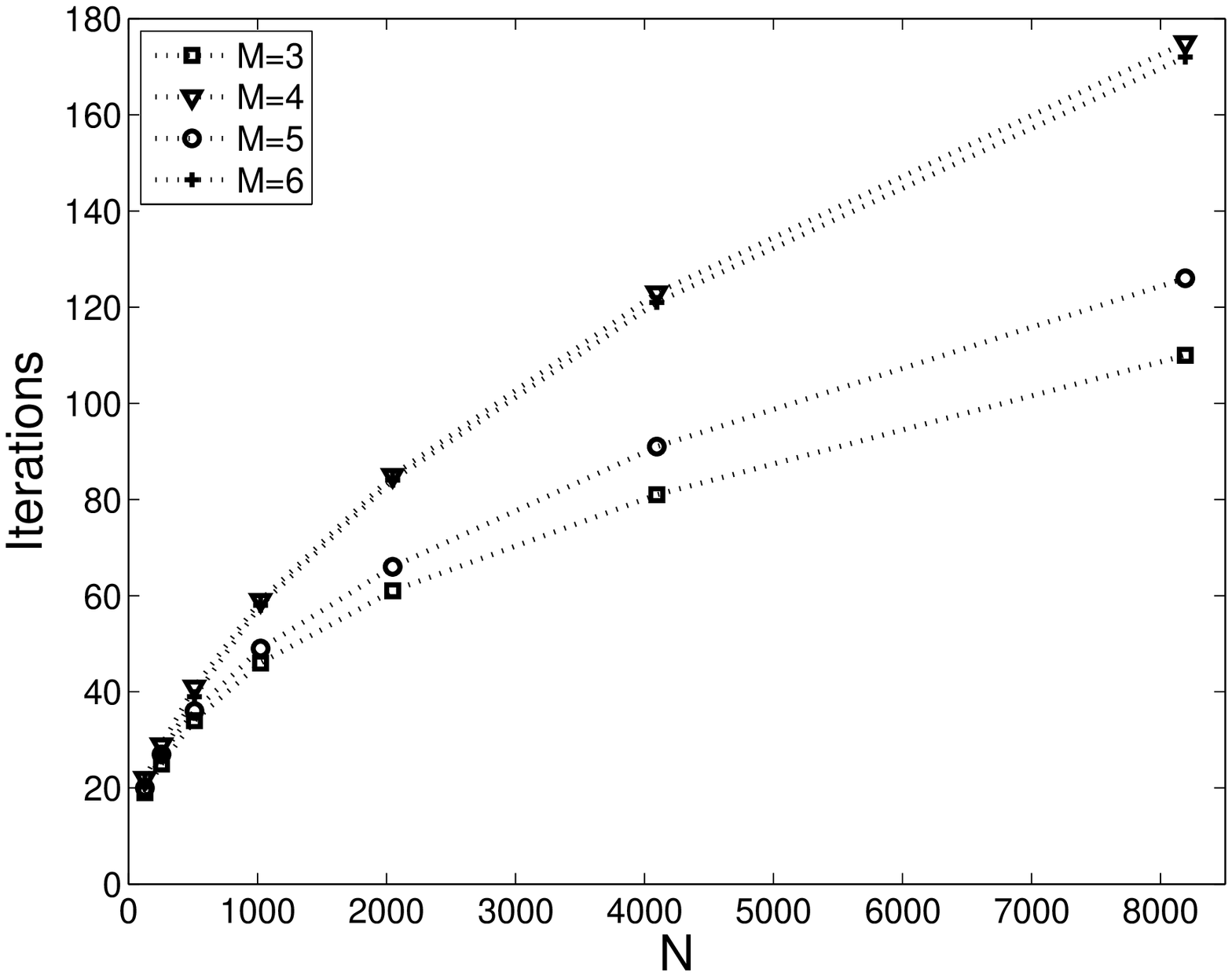}
  \caption{Convergence history ($M=3$, left) and iterations (right) of
    the Couette flow for $\Kn=0.1199$ and $u^W=1.2577$ on non-uniform
    grids.}
  \label{fig:couette-Kn01-uw300-asinh-mesh-res}
\end{figure}

Then we test the same case on a sequence of non-uniform meshes
generated by the inverse hyperbolic sine function as
\begin{align}
  \label{eq:asinh-mesh}
  x_i = x_0 + \frac{\sinh^{-1}(-5 + \frac{10 i}{N}) - \sinh^{-1}(-5)}{
    \sinh^{-1}(5) - \sinh^{-1}(-5)}, \quad i=0,1,\ldots,N.
\end{align}
This mesh is a particular setup for resolving the Knudsen layer around
the boundary. As can be seen in Table \ref{tab:couette-asinh-mesh} and
Figure \ref{fig:couette-Kn01-uw300-asinh-mesh-res}, the NMG solver
behaves similar features on these meshes. It is clear that these are
meshes seriously deviating from a quasi-uniform mesh, indicating the
NMG solver is quite stable.

\begin{table}[!tb]
  \centering
  {\begin{tabular}{l|l|ccccccc}
 \hline
\multicolumn{2}{c|}{$N$} & $2^7$ & $2^8$ & $2^9$ & $2^{10}$ & $2^{11}$ & $2^{12}$ & $2^{13}$ \\
\hline \multirow{4}{*}{\small Iterations}
& $M=3$ & 19& 25& 34& 46& 61& 81& 110 \\
& $M=4$ & 22& 29& 41& 59& 85& 123& 175\\
& $M=5$ & 20& 27& 36& 49& 66& 91& 126 \\
& $M=6$ & 21& 27& 39& 58& 84& 121& 172\\
\hline
  \end{tabular}}
\caption{Iterations of the Couette flow for $\Kn = 0.1199$ and $u^W = 1.2577$ on non-uniform grids.}
  \label{tab:couette-asinh-mesh}
\end{table}

At last, we examine the performance of the NMG solver for different
Knudsen numbers and $u^W$. Again similar convergence histories are
obtained, and the total NMG iterations are presented in Table
\ref{tab:couette-various-Kn-uw}. As the time-stepping scheme, the NMG
solver converges slower in both cases for $\Kn = 0.01199$ and $\Kn =
1.199$ than the case for $\Kn = 0.1199$. However, there is still a
substantial gain in efficiency in comparison to the single grid
solver. For large plate velocity of $u^W = 4.1923$, the NMG solver
converges a little slower than the case $u^W = 1.2577$.

\begin{table}[!tb]
  \centering
  {\begin{tabular}{l|c|c|c}
 \hline
$u^W$ & \multicolumn{2}{c|}{$1.2577$} & $4.1923$ \\
\cline{1-4}
$\Kn$ & $0.01199$ & $1.199$ & $0.1199$  \\
\hline
Iterations & 115& 120& 52 \\
\hline
  \end{tabular}}
\caption{Iterations of the Couette flow for various Knudsen numbers and $u^W$ on a uniform grid with $N=2048$ and $M=3$.}
  \label{tab:couette-various-Kn-uw}
\end{table}

\subsection{Force driven Poiseuille flow}
\label{sec:num-ex-poiseuille}
\begin{figure}[!tb]
{  \centering
  \subfigure[Density, $\rho$]{\includegraphics[width=0.5\textwidth]{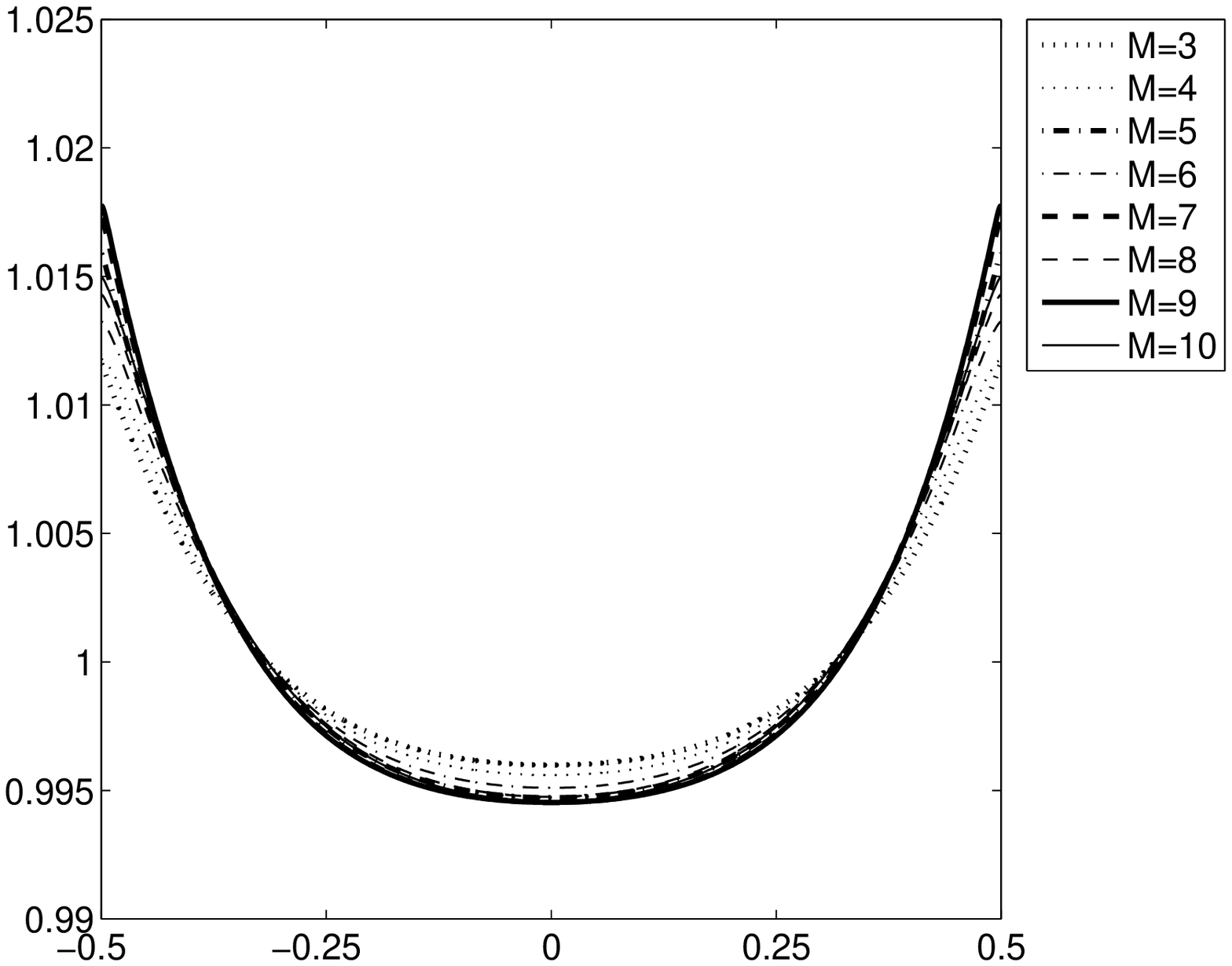}}\hfill
 \subfigure[Temperature, $\theta$]{\includegraphics[width=0.495\textwidth]{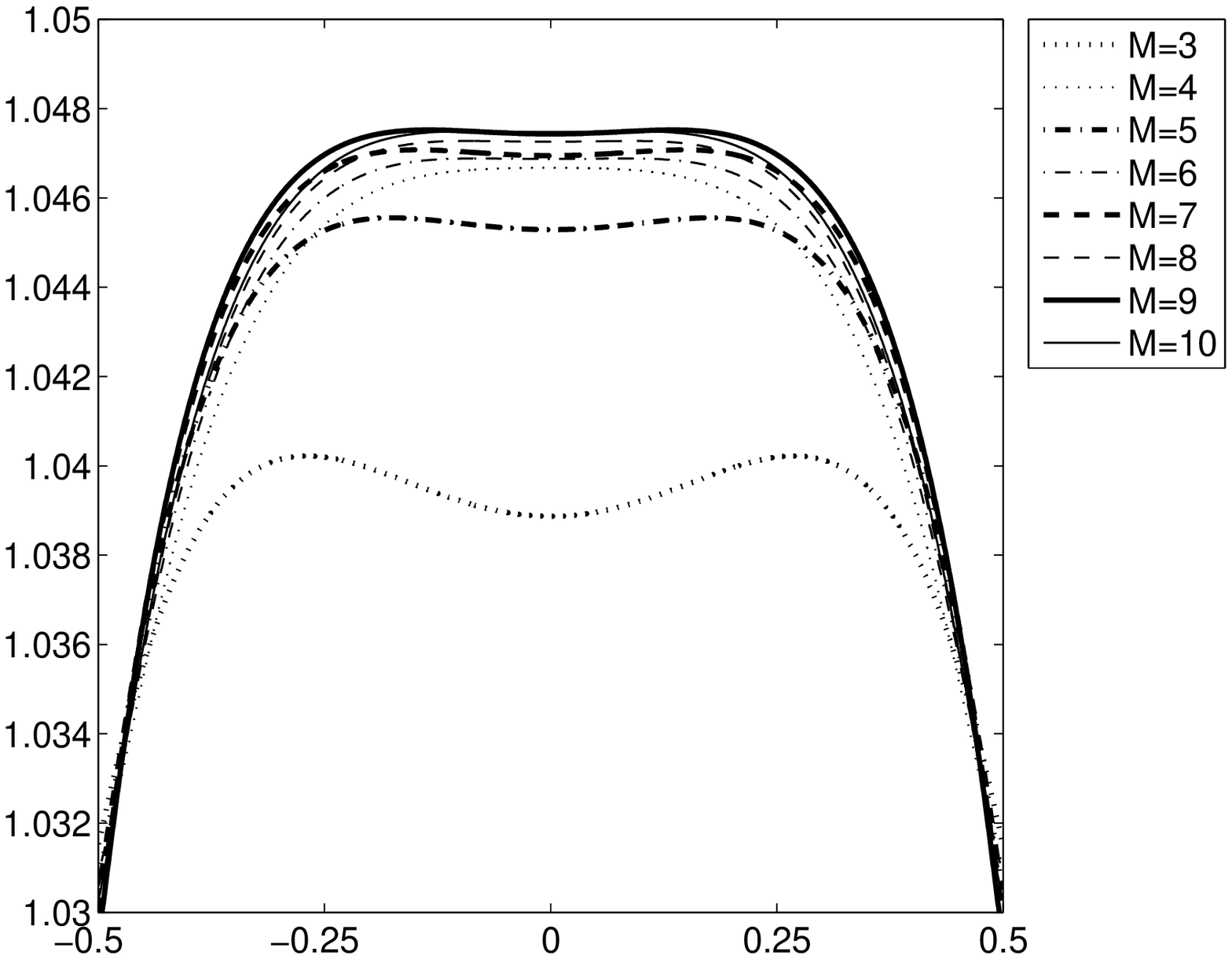}}\\
\subfigure[Normal stress, $\sigma_{11}$]{ \includegraphics[width=0.495\textwidth]{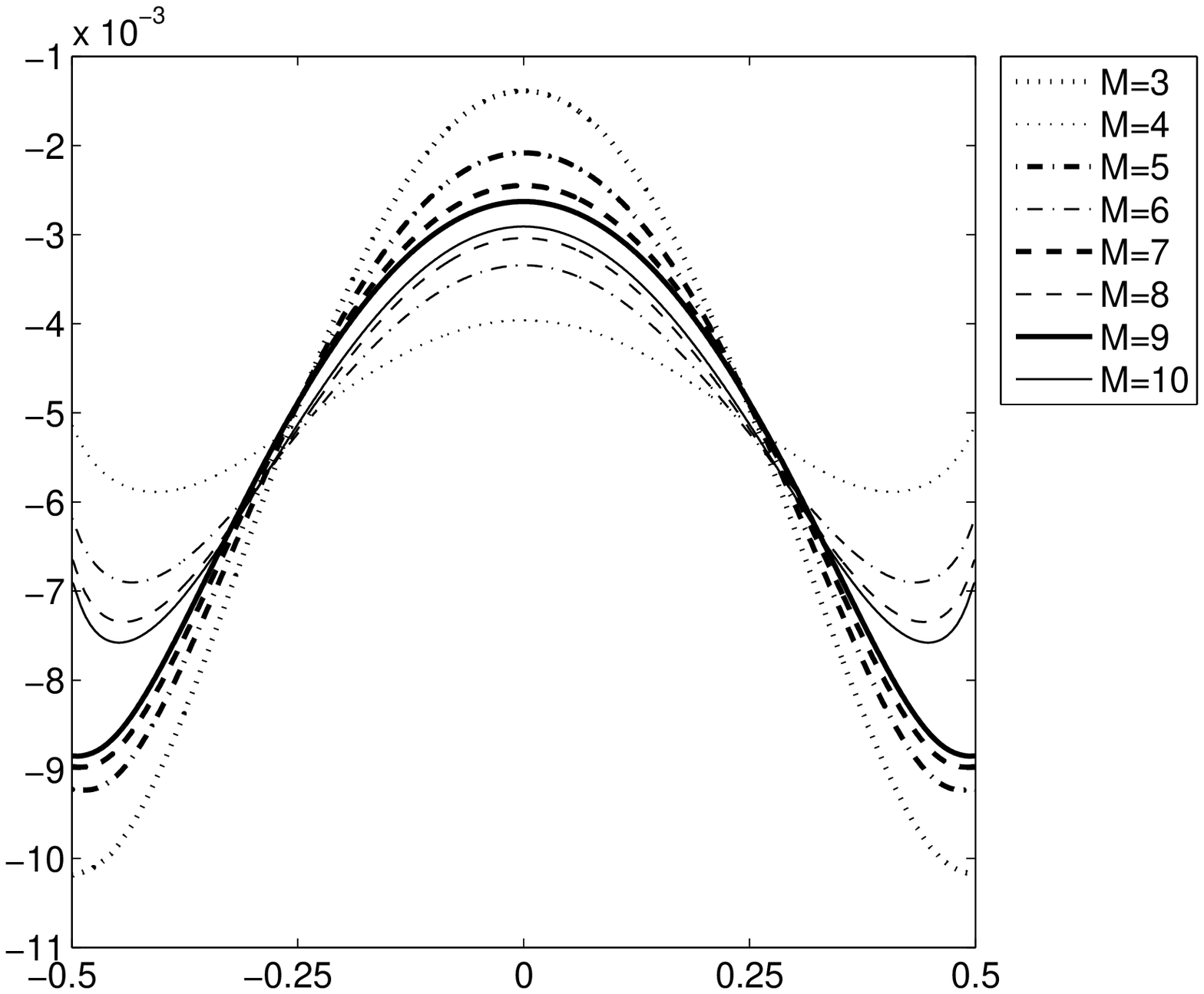}}\hfill
 \subfigure[Heat flux, $q_2$]{\includegraphics[width=0.495\textwidth]{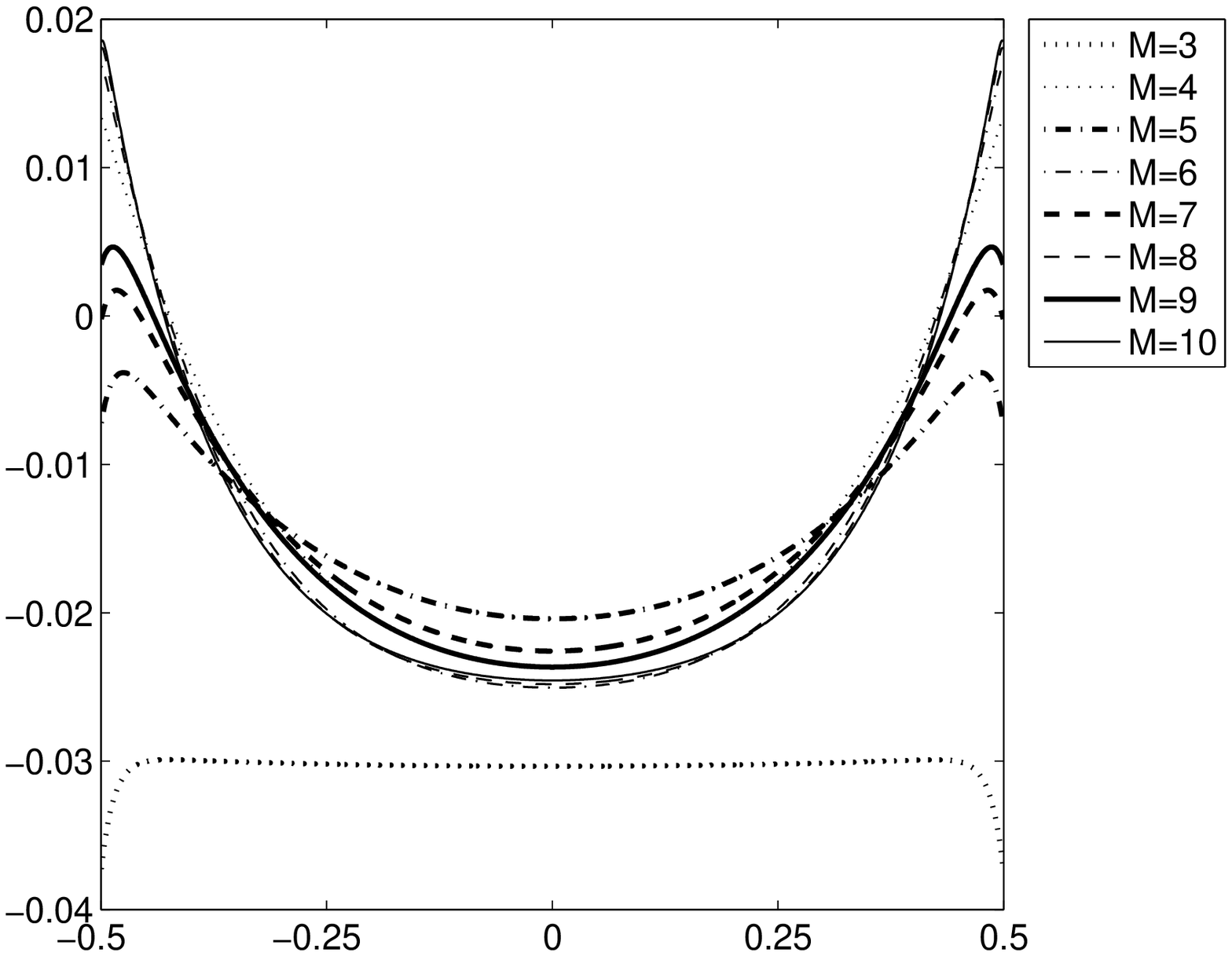}} }
  \caption{Solution of the force driven Poiseuille flow.}
  \label{fig:poiseuille-density}
\end{figure}

Next we consider the force driven Poiseuille flow which is also
frequently investigated in the literatures \cite{Garcia, Li,
  Xu2007}. For this example, the gas lies between two stationary
plates parallel to the $yz$-plane with a distance of $L=1$, and two
plates have the same temperature of $\theta^W=1$. In contrast to the
Couette flow, the Poiseuille flow is driven by an external constant
force, which is set as $\bF=(0, 0.2555,0)^T$ in our
tests. Additionally, the collision frequency $\nu$ is given by the
hard sphere model as
\begin{align}
  \label{eq:poiseuille-nu}
  \nu = \frac{16}{5}\sqrt{\frac{\theta}{2 \pi}} \frac{\Pr}{\Kn} \rho,
\end{align}
and the Knudsen number $\Kn = 0.1$ is considered. The computations
also begin with the global equilibrium \eqref{eq:couette-initial} as
the Couette flow. For these settings of the Poiseuille flow, the
solution of the Boltzmann equation using the direct simulation of
Monte Carlo (DSMC) was investigated in \cite{Garcia}, and the solution
of the Boltzmann equation with the Shakhov collision term using the
\NRxx method was considered in \cite{Li}. Here we present the solution
of the hyperbolic moment system for the Boltzmann equation with the
ES-BGK collision term in Figure \ref{fig:poiseuille-density}.

\begin{figure}[!tb]
  \centering
  \includegraphics[width=0.8\textwidth]{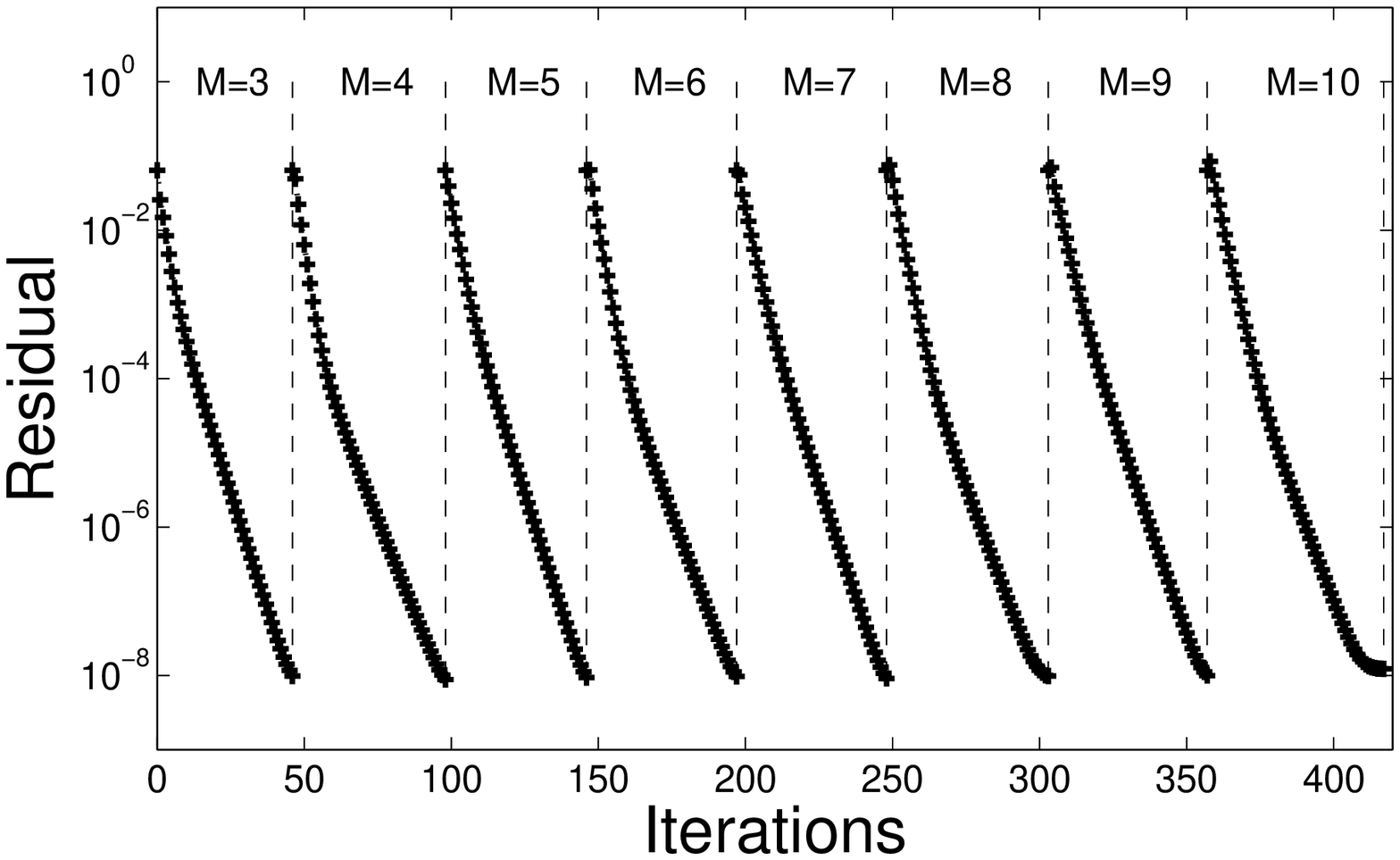}
  \caption{Convergence history of the force driven Poiseuille flow
    on the uniform grid with $N = 2048$.  }
  \label{fig:poiseuille-order-res}
\end{figure}

\begin{figure}[!tb]
  \centering
  \includegraphics[width=0.5\textwidth]{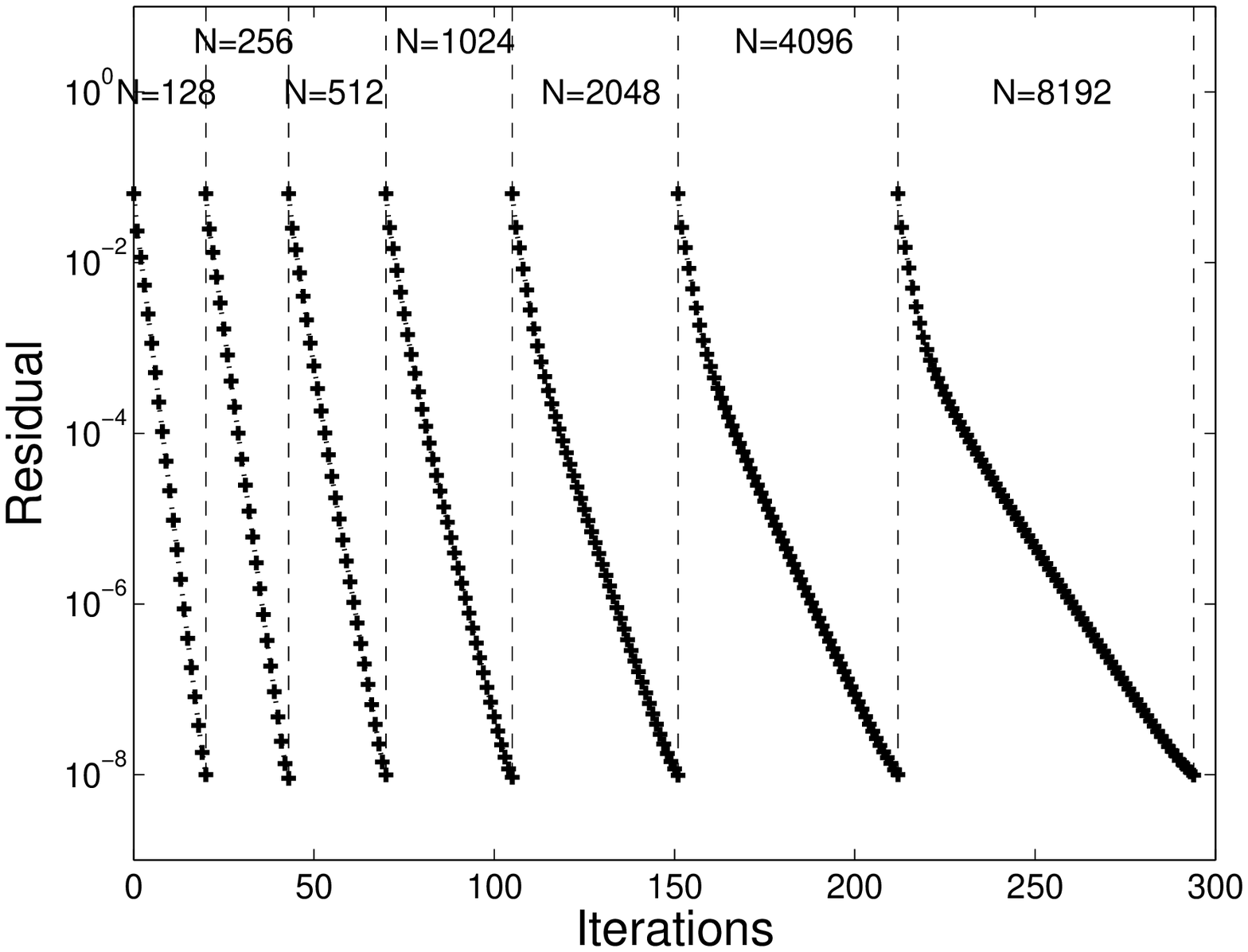}\hfill
  \includegraphics[width=0.5\textwidth]{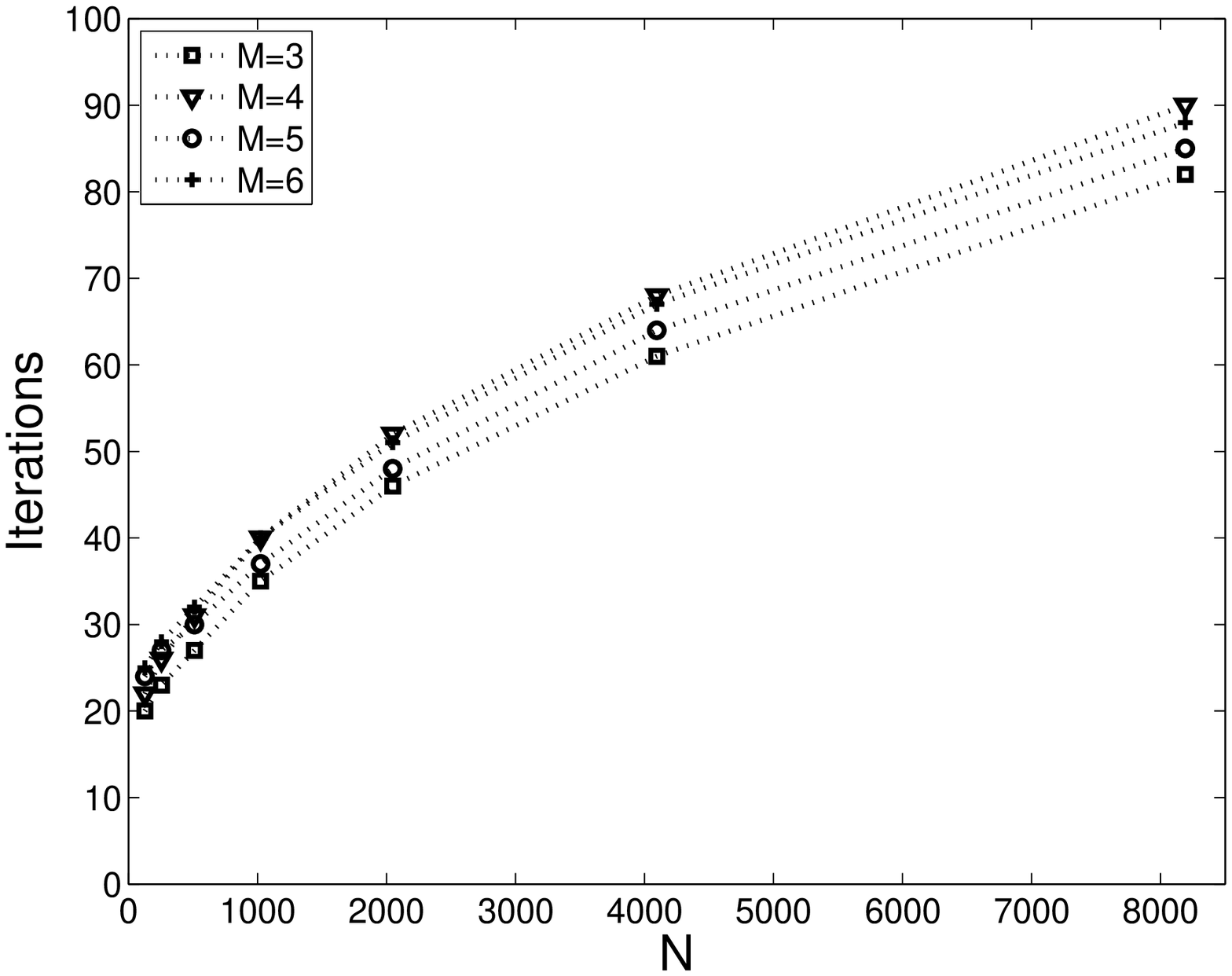}
  \caption{Convergence history ($M=3$, left) and iterations (right) of
    the force driven Poiseuille flow on a uniform grids series.  }
  \label{fig:poiseuille-mesh-res}
\end{figure}

We still omit the discussion on the solution in comparison with the
solution given in \cite{Li, Garcia}, and focus on the efficiency of
the NMG solver. The convergence histories for various $M$ on a uniform
grid with $N=2048$ and for $M=3$ on a uniform refined grids series are
shown in Figure \ref{fig:poiseuille-order-res} and Figure
\ref{fig:poiseuille-mesh-res} (left), respectively, while the
corresponding number of NMG iterations can be found in Table
\ref{tab:poiseuille-mesh}. The total iterations in terms of grid size
are presented in Figure \ref{fig:poiseuille-mesh-res} (right). As
expected, the results show the similar convergence rates for all
simulations as the force-free Couette flow, which implies a great
improvement in efficiency compared with the single grid solver.

It is also seen in Figure \ref{fig:poiseuille-order-res} that the
convergence rate slows down quickly in the last few iterations for
$M=10$. Similar situations occur in simulations of the Couette flow
with large $M$ or large $u^W$. The reason might be that in these cases
the amount of residual reduced by the smoothing steps is in the same
order as the error introduced by transferring the solution between two
successive grids. Just increasing the smoothing steps can preserve the
convergence rate during the total iterations. However, the use of
large smoothing steps during the total iterations is not preferable in
considering the computational cost. In fact, the case of the smoothing
steps $\nu_1=\nu_2=1$ is more efficient than the case of $\nu_1,\nu_2
> 1$ if it can preserve the convergence rate. Consequently, an
adaptive choice of the smoothing steps might be considered to save the
computational cost while maintaining the convergence rate.

\begin{table}[!tb]
  \centering
  {\begin{tabular}{l|l|ccccccc}
 \hline
\multicolumn{2}{c|}{$N$} & $2^7$ & $2^8$ & $2^9$ & $2^{10}$ & $2^{11}$ & $2^{12}$ & $2^{13}$ \\
\hline \multirow{4}{*}{\small Iterations}
& $M=3$ & 20 & 23 & 27 & 35 & 46 & 61 & 82\\
& $M=4$ & 22 & 26 & 31 & 40 & 52 & 68 & 90 \\
& $M=5$ & 24 & 27 & 30 & 37 & 48 & 64 & 85 \\
& $M=6$ & 25 & 28 & 32 & 40 & 51 & 67 & 88 \\
\hline
  \end{tabular}}
\caption{Iterations of the force driven Poiseuille flow on a uniform grids series.}
  \label{tab:poiseuille-mesh}
\end{table}


\section{Concluding remarks}
\label{sec:conclusion}
An efficient and robust nonlinear multigrid solver has been developed
for the hyperbolic moment system derived for the steady-state
Boltzmann equation with ES-BGK collision term. The moment system is
discretized using the unified framework of the \NRxx method such that
our solver is also unified for the system with moments up to arbitrary
order. A nonlinear iterative method, namely, SGS-Newton iteration, is
presented to solve the resulting discretized system on a single grid
level. It is then accelerated drastically by putting itself as
smoother into a routine nonlinear multigrid procedure. The numerical
experiments on two benchmark problems demonstrate the efficiency and
the robustness of our NMG solver.

Although we have considered only the moment system for the Boltzmann
equation with ES-BGK collision term, the implementation of the
proposed NMG solver is clearly not restricted to this collision
model. It is trivial to extend our NMG solver to the BGK model since
it is a special case of the ES-BGK model with $\Pr=1$.  For the
Shakhov model, the extension is still quite simple with the expansion
of the collision term given in \cite{Li}.

The extension of the NMG solver to high spatial dimensional case is
under our current study and will be reported elsewhere.

\section*{Acknowledgements}

The research of Z. Hu was supported in part by the China Postdoctoral
Science Foundation (2013M540807). The research of R. Li was supported
in part by the National Basic Research Program of China (2011CB309704)
and the National Science Foundation of China (11325102, 91330205).


\appendix
\section*{Appendices}

\section{Computation of the parameter $\hat{\tau}$}
\label{sec:def-tau}
Suppose the distribution function $g$ belongs to the space
$\mF_M(\bu,\theta)$, where the relation \eqref{eq:moments-relation}
holds for its coefficients $g_{\alpha}$. When approximate $g$ in
another space $\mF_M(\bu',\theta')$, the corresponding coefficients,
denoted by $g_{\alpha}'$, would not satisfy
\eqref{eq:moments-relation} usually. However, similar relation can be
deduced from \eqref{eq:moments}, which is given as
\begin{align}
  \label{eq:moments-relation-other-basis}
  \begin{aligned}
    & \rho = g_0', \\
    & \rho ( \bu - \bu' ) = (g_{e_1}', g_{e_2}', g_{e_3}')^T, \\
    & \rho |\bu - \bu' |^2 + 3 \rho (\theta - \theta') = 2
    \sum_{d=1}^3 g_{2e_d}'.    
  \end{aligned}
\end{align}
The last two equations of \eqref{eq:moments-relation-other-basis}
indicate
\begin{align}
  \label{eq:new-theta-relation}
  \theta = \frac{2\sum_{d=1}^3 g_{2e_d}' - \sum_{d=1}^3 (g_{e_d}' )^2
    / \rho}{3 \rho} + \theta'.
\end{align}

Now corresponding to the step of updating solution by
\eqref{eq:update-dis} in the local Newton iteration (Algorithm
\ref{alg:local-newton}), we have
\begin{align}
 \label{eq:relations-update-dis}
\begin{aligned}
& \rho = g_0' = f_{i,0}^{(m)} + \tau \Delta f_{i,0}^{(m)}, \qquad
 \bu' = \bu_i^{(m)}, \qquad \theta' = \theta_i^{(m)},\\
& g_{e_d}' = \tau \Delta f_{i,e_d}^{(m)}, \quad d=1,2,3, \qquad
 \sum_{d=1}^3 g_{2e_d}' = \tau \sum_{d=1}^3 \Delta f_{i,2e_d}^{(m)}.
\end{aligned}
\end{align}
Substituting the above equations into \eqref{eq:new-theta-relation}
immediately yields
\begin{align}
  \label{eq:new-theta}
  \theta = \frac{ 2\tau \sum_{d=1}^3 \Delta f_{i,2e_d}^{(m)} - \tau^2
    \sum_{d=1}^3 (\Delta f_{i,e_d}^{(m)} )^2 / (f_{i,0}^{(m)} + \tau
    \Delta f_{i,0}^{(m)})}{3 (f_{i,0}^{(m)} + \tau \Delta
    f_{i,0}^{(m)})} + \theta_i^{(m)}.
\end{align}
In our implementation, the positivity of the density and the
temperature are preserved during the iterations, that is,
\begin{align*}
  \rho \geq \bar{\rho} >0, \qquad  \theta \geq \bar{\theta} >0,
\end{align*}
where $\bar{\rho}$, $\bar{\theta}$ are the given lower bounds of the
density and the temperature respectively. 

For the positivity of the density, we have
\begin{align}
  \label{eq:tau-density}
0<  \tau \leq \frac{\bar{\rho} - f_{i,0}^{(m)}}{ \Delta
    f_{i,0}^{(m)}},
\end{align}
if $\Delta f_{i,0}^{(m)} < 0$. Otherwise, the density is always
positive for $\tau>0$.

For the positivity of the temperature, we deduce from
\eqref{eq:new-theta} that
\begin{align}
  \label{eq:tau-theta-inequality}
A \tau^2 + B \tau + C \leq 0, 
\end{align}
where 
\begin{align*}
  & A = \sum_{d=1}^3 \left( \Delta f_{i,e_d}^{(m)}\right)^2 - 2 \Delta
  f_{i,0}^{(m)} \sum_{d=1}^3 \Delta f_{i,2e_d}^{(m)} - (
  \theta_i^{(m)} - \bar{\theta} ) \left(\Delta f_{i,0}^{(m)}\right)^2,
  \\ & B = -2 f_{i,0}^{(m)} \left(\sum_{d=1}^3 \Delta f_{2e_d}^{(m)} +
    \Delta f_{i,0}^{(m)} (\theta_i^{(m)} - \bar{\theta} )\right), \qquad
  C = -(f_{i,0}^{(m)})^2 (\theta_i^{(m)} - \bar{\theta}).
\end{align*}
It is trivial to solve the above inequality, and the solution is given
as
\begin{enumerate}[i).]
\item $ 0 < \tau \leq -\frac{C}{B}$ if $A=0$, $B>0$.
\item $ 0 < \tau \leq \frac{-B + \sqrt{B^2-4A C}}{2A}$, if $A> 0$,
  $B^2-4A C > 0$ or $A<0$, $B>0$, $B^2-4AC>0$.
\item Otherwise, $\tau > 0$.
\end{enumerate}

Finally, the parameter $\hat{\tau}$ is obtained such that the above
inequalities \eqref{eq:tau-density} and \eqref{eq:tau-theta-inequality}
hold for $\tau$.

\section{Construction of the restriction operator $I_h^H$}
\label{sec:def-restriction}
For any fine grid function $g_{\sss h}$ with $g_{\sss h,i} \in
\mF_M(\bu_{\sss h,i},\theta_{h,i})$, denote its restriction $I_h^H
g_{\sss h}$ by $g_{\sss H}$. Obviously, it is enough to construct
$g_{\sss H,i}$ on the $i$-th coarse grid cell $[x_{\sss H,i},x_{\sss
  H,i+1}]$. Suppose $g_{\sss H,i}$ belongs to $\mF_M(\bu_{\sss
  H,i},\theta_{H,i})$, where $\bu_{\sss H,i}$ and $\theta_{H,i}$ are
macroscopic velocity and temperature of the restriction $I_h^H f_h$ on
the $i$-th coarse grid cell, in which $f_h$ is the fine grid solution
for \eqref{eq:multi-fine-problem}.

Due to the importance of conservation, $g_{\sss H,i}$ is required to
preserve this property. To be specific, the following equation
\begin{align}
  \label{eq:conservative-restriction}
  & \int_{x_{H,i}}^{x_{H,i+1}} \int g_{\sss H, i}(\bxi) p(\bxi) \dd
  \bxi \dd x = \int_{x_{h,2i}}^{x_{h,2i+1}} \int g_{\sss h, 2i}(\bxi)
  p(\bxi) \dd \bxi \dd x + \int_{x_{h,2i+1}}^{x_{h,2i+2}} \int g_{\sss
    h, 2i+1}(\bxi) p(\bxi) \dd \bxi \dd x
\end{align}
should hold for any polynomial $p(\bxi)$ of degree no more than
$M$. To evaluate an arbitrary restriction $g_{\sss H,i}$, one should
first calculate $\bu_{\sss H,i}$ and $\theta_{H,i}$. To this end,
replacing $g_{\sss H}$ and $g_{\sss h}$ by the solution $f_{H}$ and
$f_{h}$ respectively, and employing \eqref{eq:moments}, we have
\begin{align}
  \label{eq:restriction-rho-u-theta}
\begin{aligned}
  \rho_{\sss H,i} \Delta x_{\sss H,i} = &~ \rho_{\sss h,2i} \Delta
  x_{\sss h,2i} + \rho_{\sss h,2i+1} \Delta x_{\sss h,2i+1},\\
  \rho_{\sss H,i} \bu_{\sss H,i} \Delta x_{\sss H,i} = &~ \rho_{\sss
    h,2i} \bu_{\sss h,2i} \Delta x_{\sss h,2i} + \rho_{\sss h,2i+1}
  \bu_{\sss h,2i+1} \Delta x_{\sss h,2i+1},\\ \left( \rho_{\sss H,i}
    \bu_{\sss H,i}^2 + 3\rho_{\sss H,i}\theta_{H,i} \right) \Delta
  x_{\sss H,i} = &~ \Delta x_{\sss h,2i} \left( \rho_{\sss h,2i}
    \bu_{\sss h,2i}^2 + 3 \rho_{\sss h,2i}\theta_{h,2i} \right) \\ & +
  \Delta x_{\sss h,2i+1} \left( \rho_{\sss h,2i+1} \bu_{\sss h,2i+1}^2
    + 3 \rho_{\sss h,2i+1}\theta_{h,2i+1} \right).
\end{aligned}
\end{align}
After computing $\bu_{\sss H,i}$, $\theta_{H,i}$ from the above equations,
we employ the transformation \cite{NRxx} to project $g_{\sss h, 2i}$
and $g_{\sss h, 2i+1}$ into $\mF_M(\bu_{\sss H,i},\theta_{H,i})$, denoted
by $\tilde{g}_{\sss h, 2i}$ and $\tilde{g}_{\sss h, 2i+1}$
respectively. As the transformation is conservative,
\eqref{eq:conservative-restriction} is re-written as
\begin{align}
  \label{eq:conservative-restriction-same-basis}\begin{aligned}
    \int_{x_{H,i}}^{x_{H,i+1}} \int g_{\sss H, i}(\bxi) p(\bxi) \dd
    \bxi \dd x = & \int_{x_{h,2i}}^{x_{h,2i+1}} \int \tilde{g}_{\sss
      h, 2i}(\bxi) p(\bxi) \dd \bxi \dd x \\ & +
    \int_{x_{h,2i+1}}^{x_{h,2i+2}} \int \tilde{g}_{\sss h, 2i+1}(\bxi)
    p(\bxi) \dd \bxi \dd x.
  \end{aligned}
\end{align}

Now let
\begin{align}
  \label{eq:test-basis}
  p_\alpha(\bxi)=\mH_{\theta_{H,i}, \alpha}\left( \frac{\bxi -
      \bu_{\sss H,i}}{\sqrt{ \theta_{H,i} }} \right) \exp\left(
    \frac{|\bxi - \bu_{\sss H,i}|^2}{2 \theta_{H,i}} \right), \quad
  |\alpha| \leq M.
\end{align}
It is trivial to show that $\{p_\alpha(\bxi)\}_{|\alpha|\leq M}$ forms
a complete orthogonal basis of the polynomial space of degree no more
than $M$. Substituting \eqref{eq:test-basis} into
\eqref{eq:conservative-restriction-same-basis}, and employing the
orthogonality of the basis, we deduce
\begin{align}
  \label{eq:restriction-moments-conservation}
  \int_{x_{H,i}}^{x_{H,i+1}} g_{\sss H,i,\alpha} \dd x =
  \int_{x_{h,2i}}^{x_{h,2i+1}} \tilde{g}_{\sss h,2i,\alpha} \dd x +
  \int_{x_{h,2i+1}}^{x_{h,2i+2}} \tilde{g}_{\sss h,2i+1,\alpha} \dd x,
  \quad |\alpha|\leq M,
\end{align}
which follows that
\begin{align}
  \label{eq:restriction-moments}
  g_{\sss H,i,\alpha} = \frac{\tilde{g}_{\sss h,2i,\alpha}\Delta
    x_{\sss h,2i} + \tilde{g}_{\sss h,2i+1,\alpha} \Delta x_{\sss
      h,2i+1}}{\Delta x_{\sss H,i}}, \quad |\alpha| \leq M.
\end{align}


\bibliographystyle{plain}
\bibliography{../article,../tiao}
\end{document}